\documentclass{article}
\usepackage{amsmath,amssymb,amsthm,bm}
\usepackage{caption,graphicx}
\usepackage[all]{xy}
\usepackage{hyperref}
\theoremstyle{plain}
\newtheorem{thm}{Theorem}
\newtheorem{lm}{Lemma}
\newtheorem{prop}{Proposition}
\newtheorem{cor}{Corollary}
\def\be{\begin{equation}}
\def\ee{\end{equation}}

\title{On the $abc$ conjecture}
\date { \ }
\author{C. {\sc Soul\'e}\footnote{CNRS and  IH\'ES, Le Bois-Marie, 35 route de Chartres, 91440 {\sc Bures-sur-Yvette}, France \newline soule@ihes.fr}}

\begin{document}

\maketitle

\bigskip

%\hfill     Pour Olivier

\bigskip

{\centerline{\bf Introduction}}

\medskip

Let $K_0$ be a number field and ${\mathcal C}$ a smooth, projective, geometrically irreducible curve of genus $2$ over $K_0$.  We denote by
$$
h : {\mathcal C} (\overline{K_0}) \to {\mathbb R}
$$
a logarithmic height associated to the canonical line bundle of ${\mathcal C}$. In this paper we shall prove the following effective version of the  Mordell conjecture for ${\mathcal C}$, a weak form of Vojta's conjecture \cite{V}, Conj. 5.2.6:

\bigskip
\noindent {\bf Theorem A.} {\it Given any positive real number $\varepsilon > 0$, there exists a constant $\beta (\varepsilon) > 0$ such that, for any number field $F$ containing $K_0$, of degree $[F:{\mathbb Q}]$ and discriminant $\Delta_F$ over ${\mathbb Q}$, and for any rational point $P \in {\mathcal C}(F)$, the following inequality holds:
$$
h(P) \leq \left( 12 + \varepsilon \right)  \frac{\log \vert \Delta_F \vert}{[F:{\mathbb Q}]} + \beta (\varepsilon) \, .
$$
}

\bigskip

The $abc$ conjecture of Masser and Oesterl\'e is known to follow from Theorem~A  \cite{MB3}, albeit with a wrong exponent
(compare \cite{Oe}  Conj. 3):

\bigskip

\noindent {\bf Theorem B.} {\it  Given any positive real number $\varepsilon > 0$, there exists a positive real number $A(\varepsilon) $ such that, if $a,b,c$ are positive integers with no common divisor and such that
\be
\nonumber
a+b=c \, ,
\ee
then
\be
\nonumber
c \leq A(\varepsilon) \, N^{24 + \varepsilon } \, ,
\ee
where $N$ is the product of the prime numbers dividing $abc$.
}

\medskip

\medskip

\newpage

Our proof of Theorem~A follows the strategy proposed by Parshin \cite{Pa1} and then Moret-Bailly \cite{MB3}. 
Namely, a construction of Kodaira \cite{K} (see also \cite{MD}) gives, after choosing a finite extension $K$ of $K_0$, 
 a finite \'etale map $Y \to {\mathcal C} \underset{K_0}{\otimes} K$ and a non-isotrivial smooth proper morphism
$
\phi : X \to Y
$
with fibers geometrically connected curves of genus $6$. If $F$ is a finite extension of $K$ and if $P \in Y (F)$ is a rational point, we let $X_{P}$ be the fiber of $\phi$ at $P$. Given two  line bundles $\overline M_1$ and $\overline M_2$ on $X_P$, equipped with adelic metrics (see  \cite{Z2}, \cite {YZ}, and \S 1.3.2), one can define an arithmetic  intersection number $\overline{M_1} \cdot \overline{M_2} \in {\mathbb R}$. This pairing extends the Arakelov intersection theory on arithmetic surfaces \cite{A}. Let $\omega$ be the sheaf of relative differential forms of $X$ over $Y$. We endow $  \omega$ with an adelic metric and,  for each point $P$ as above,  we let
$\overline {\omega_P}$ be the restriction of $\overline {\omega}$ to $X_P$. Parshin and Moret-Bailly noticed that the 
map associating to $P$ the real number $\overline {\omega_P}^2=\overline {\omega_P}\cdot\overline {\omega_P}$
 is a logarithmic height  on $Y$. Therefore, to prove Theorem~A it is enough to  relate heights on $\mathcal C$ to heights on $Y$
 and to bound $\overline\omega_P^2$ from above.

\smallskip

To bound $\overline\omega_P^2$ we proceed as follows. First, we
  choose a  line bundle $L_P$  on $X_P$ that we endow with
an adelic metric. The intersection number ${\overline L_P}\cdot {\overline L_P}$
is, up to a bounded constant, a multiple of $\overline\omega_P^2$.
On the other hand, after some base change, Minkowski's theorem provides
an upper bound for  ${\overline L_P}\cdot {\overline L_P}$. The combination
of these two estimates of  ${\overline L_P}\cdot {\overline L_P}$ leads to an upper bound for
$\overline\omega_P^2$.

\smallskip

The paper is organized as follows. In Section~1 we fix notations and we prove
that the sup-norm  on the global sections of an  adelic line bundle on a curve is adelic 
( Proposition~\ref{G}). We prove also 
a general result, Proposition~\ref{prop1}, which 
will enable us to deduce inequalities for  arithmetic intersection numbers  on a curve 
(up to bounded constants) from inequalities among classical intersection numbers. Next, in Section~2, given
two ample line bundles $L_1$ and $L_2$ on a curve equipped with {\it admissible \it} metrics, we derive an upper
bound for ${\overline L_1}\cdot {\overline L_2}$ (Proposition~\ref{prop12}).
 In Section~3, with the help of Voisin, we compute Betti numbers on some curves  of genus $3$ and $6$.
 In Section~4  we consider stable 2-pointed curves
over a number field $F$ and we show, following Abramovich, that, after an alteration of the base, they admit
a model  over the integers of $F$ (Proposition~\ref{prop19}).
In Section~5 we introduce the Kodaira  fibration $\phi : X \to Y$
and we compute some classical
intersection numbers on $X$. 

\smallskip

In Section~6, we  state the main theorem,
Theorem~\ref{thm1}, which gives an upper bound for $\overline\omega_P^2$ in terms of the logarithm
of the absolute discriminant of the  field of definition of $P$. To prove that theorem
we define  a line bundle $L$  on $X$. We compute
the classical intersection number $L\cdot L$ and we bound from
below the degree of $\phi_*(L)$  (Proposition~\ref{prop6} ).
We then bound  ${\overline L_P}\cdot {\overline L_P}$ from above by first using
Proposition~\ref{prop12}, next, after the base change of Proposition~\ref{prop19},  applying  Proposition~\ref{G} and the Minkowski's
theorem on successive minima, and finally using Proposition~\ref{prop6} and  Proposition~\ref{prop1}.
We conclude by showing in § 6.3.7 that ${\overline L_P}\cdot {\overline L_P}$ is a multiple of $\overline\omega_P^2$ 
(up to a constant).

\smallskip

In Section~7, following \cite{MB3}, we deduce some  results from Theorem~\ref{thm1}. First, 
we prove Theorem~A (see Theorem~\ref{thm2} ).  Next we prove a weak version of Szpiro's conjecture of the
discriminant, Theorem~\ref{thmDI}.Then we prove Theorem~B and, more generally,
 a  weak uniform $abc$ theorem for number fields, Theorem~\ref{thm3}. The consequences
of the $abc$ conjecture  follow also from Theorem~\ref{thm3}, with weaker constants.
For instance, following \cite{G-St}, we deduce from Theorem~\ref{thm3} that the $L$-function of a quadratic character with  negative  discriminant does not have any Siegel zero  (Theorem~\ref{thm77}) (see also \cite{OA} and \cite{KE}).
\medskip

I thank D. Abramovich ,   A. Beauvillle, S. Bloch, J-B. Bost, J. Boxall, F. Catanese, A. Chambert-Loir,
 J-L. Colliot-Th\'el\`ene,   G. Faltings,  O. Gabber, C. Gasbarri,  L. Gatto,  H. Gillet,  
  W.  Gubler, Q. Guignard, J. Lannes, Q. Liu,  L. Moret-Bailly, J. Oertel,  P. Scholze, J-P. Serre,
T. Shaska,  M. Stoll,  C. Voisin, P. Vojta and X. Yuan for their help.

\bigskip

\noindent {\bf Notation :} If $A$ is an abelian group we let 
$A_{\mathbb Q} = A \underset{\mathbb Z}{\otimes} \, {\mathbb Q}$. In the following $K$ is a  number field. 
When an adelic metric has been chosen on a line bundle $M$ on some variety over $K$, we let $\overline M$ be the corresponding adelic line bundle. 
%%%%
\section{Preliminaries}

\subsection{ \\} We denote by $d_K=[K : {\mathbb Q}]$ the absolute degree of $K$ and by $\Delta_K$ its absolute discriminant. We write
$$
\delta(K) = \log (\vert \Delta_K \vert) \, .
$$

Let $F$ be a finite extension of $K$, let $Y$ be a proper variety over $K$  and let
$P \in Y (F)$ be a rational point. We let 
\be
\nonumber
u : {\rm Spec}(F) \to Y
\ee
be the map defined by $P$. Assume that, for a class of maps $u$ as above, we can define a real number $a(u)$. If 
$d = [F : {\mathbb Q}]$ we shall write
\be
\nonumber
a(u) \leq 0(d)
\ee
when there exists a real constant $A > 0$, independent from $F$ and $P$, such that
\be
\nonumber
a(u) \leq A \, d \, .
\ee

\subsection{ \ } Let $Y$ a variety over $K$. We denote by ${\rm Pic} (Y)$ the Picard group of $Y$.  If $Y$ is a  projective curve on $K$ we let
\be
\nonumber
 {\deg}_K : {\rm Pic} (Y) \to {\mathbb Z}
\ee
be the degree map. If $Y$ is smooth and if $E$ is a coherent sheaf on $Y$ we let $\det (E)$ be the  determinant of $E$ and
$
\deg_K (E) = \deg_K (\det (E)) \, .
$

If $f : X \to Y$ is a flat projective curve over a variety $Y$, given two line bundles $L$ and $M$ over $X$, we let $\langle  L,M \rangle \in {\rm Pic} (Y)$ be the Deligne pairing of $L$ and $M$ \cite{De}. When $Y$ is a projective curve we define
\be
\nonumber
L \cdot M = {\deg}_K (\langle L,M \rangle) \in {\mathbb Z} \, .
\ee
When $M=L$ we write $L^2$ instead of $L \cdot L$.

\medskip
If $L$ and $M$ are two line bundles on a scheme $X$ we denote by $L+M$ the tensor product $L\otimes_{\mathcal O_X}M$.
If $D$ is Cartier divisor on $X$ we denote by $D+L$ the line bundle $\mathcal O_X(D) \otimes _{\mathcal O_X}L$.
When $X$ is a variety on $K$, if $D$ is a Cartier divisor on $X$, $L$ a line bundle on $X$, and $n \geq 0$ an integer, we write
$$
h^n (D+L) = h^n (X,D+L) = \dim_{K} (H^n (X , {\mathcal O}_X(D)\otimes_{{\mathcal O}_X}L))\, .
$$

\subsection{\\}

\subsubsection{\\}
If  $v$ is a place of $K$ we let  $K_v$ be the completion of $K$ 
 at $v$, $\pi_v$ a uniformizer of $K_v$, and $N(v)$ the cardinality of its residue field.
For every finite extension $F$ of $K$ and every finite place $w$ of $F$ restricting to $v$ on $K$, we define a norm 
$\vert \cdot \vert_w$ on $F_{w}$ by imposing that 
$$ \vert \pi_v \vert_w = \frac {1}{N(v)} \, . $$
When $w$ is archimedean, $\vert \cdot \vert_w$ is the usual norm.
 If $\overline K_v$  is an algebraic closure of $K_v$, the norm $\vert \cdot \vert_v $ extends to $\overline K_v$,
and any $x\in K^*$ satisfies the product formula
$$\prod_v \vert x \vert_v =1 \, ,$$
where $v$ runs over all places of $K$.

\medskip

\subsubsection{\\}

 Let  $X$ be a smooth  projective variety over $K$. If  $v$ is a place of $K$ we let  $\overline K_v$  be an algebraic closure of $K_v$. An {\it adelic  line bundle}
  $\overline L = (L,\Vert \cdot \Vert)$
  over $X$  (see \cite{Z2} (1.2)) is
 an 	algebraic line bundle $L$ on $X$ endowed with a collection of norms $\Vert \cdot \Vert = \{\Vert \cdot \Vert_v\}$
  on the line bundles $L_v = L \otimes \overline K_v$  defined by $L$ on $X_v = X\otimes \overline K_v$ for every place $v$ of $K$:
for every point $x\in X(\overline K_v)$, the fiber $L_v(x)$ of $L_v$ at $x$ is equipped with a $\overline K_v$-norm, i.e., if $\lambda \in \overline{ K_v}$
and if $s$ is a section in $L_v(x)$,  one  has
$$\Vert\lambda s\Vert_v = \vert \lambda \vert_v \Vert s \Vert_v \, .$$
The norm $\Vert \cdot \Vert_v$ is assumed to be ${\rm Gal}(\overline{ K_v}/K_v)$-invariant and  m-continuous 
in the sense of \cite{YZ} A.5.1.

 \smallskip
 One assumes that the collection $\Vert \cdot \Vert = \{\Vert \cdot \Vert_v\}$ is $ coherent$:
 there exists a non empty open subset $U \subset {\rm Spec}({\mathcal O}_K)$, a projective model $X_U$
 of $X$ over $U$ and a line bundle ${ L_U}$ on $X_U$  extending  $L$ on $X$ and such
 that, for every  place $v\in U$, the metric  $\Vert \cdot \Vert_v$
is induced  by the  line bundle ${L_U} \otimes{\overline  {\mathcal O}_v}$ on ${X_U}\times_U {\rm Spec}  ({\overline {\mathcal O}_v})$,
where ${\overline {\mathcal O}_v}$ is the ring of integers in ${\overline K_v}$.

 \medskip

 Let $\widehat{{\cal P}{\rm ic}}(X)$  be the Picard category of adelic line bundles over $X$.
 If $\overline L$ and $\overline M$ are two   metrized line bundles on $X$ we denote by $\overline L+\overline M$ their tensor product.
 We let $\widehat{{\rm Pic}}(X)$  be the Picard group of isomorphim classes of adelic  line bundles on $X$.
 
 \medskip
 
 An adelic  line bundle on $X$ is called  $integrable$ if it is the quotient of two semi-positive  adelic  line bundles, see \cite{Z2} (1.5)  and \cite{YZ} 2.1.1.
\medskip

\subsubsection{\\}

 When $X = {\rm Spec} ({K})$ we define as follows the 
 {\it arithmetic degree} morphism 
 $$\widehat{\deg}_K : \widehat{\rm Pic} ({\rm Spec} ({K})) \to {\mathbb R} \, .$$ 
Let $\overline L$ be an adelic  line bundle
 on ${\rm Spec} ({K})$ and $s$ a non-trivial  element in $L$.
Then
 $$\widehat{\deg}_K(\overline L) = - \sum_v \log(\Vert s \Vert_v)    \, .$$
 The sum is finite and, by the product formula for the valuations  $|\cdot|_v$, the number $\widehat{\deg}_K (\overline L)$ is independent of the choice of $s$.

\smallskip

If $X$ is a smooth, irreducible, projective  curve over $K$ and if $\overline L$ and $\overline M$ are two 
adelic   line bundles with integrable metrics on $X$, we let $\langle \overline L , \overline M \rangle \in \widehat{{\cal P}{\rm ic}} ({\rm Spec} ({K})) $
   be the Deligne pairing of $L$ and $M$, equipped with the adelic metric defined in \cite{De} (6.3.2)
and in \cite{YZ} Thm~4.1.3. The {\it arithmetic intersection number} of $\overline L$ with $\overline M$ is
$$
\overline L \cdot \overline M = \widehat\deg_K (\langle \overline L , \overline M \rangle) \, .
$$
When $\overline M = \overline L$ we write $\overline L^2 = \overline L \cdot \overline L$.

\bigskip

\subsection{\\}
Let $X$ be a smooth projective variety over $K$. Let $X^{\rm an}$ be the Berkovich space associated
to $X$ (\cite{YZ},3.1.1) and $X^{\rm r-an} \subset X^{\rm an}$ the restricted analytic space associated
to the scheme $X/{{\mathcal O}_K}$ (\cite{YZ},3.5.1).The Picard category $\widehat{{\cal P}{\rm ic}}(X)$
of adelic line bundles on $X$ is equivalent to the category of line bundles $L$ on $X$ equipped with continuous, coherent metrics on the fibers of $L$ over $X^{\rm r-an}$ (\cite{YZ}, Theorem 3.5.2).

\medskip

\subsection{\\}

\begin{lm}\label{lm1} {\it Let $X$ be a smooth projective variety over $K$, let $L$ be a line bundle
over $X$, and let $\Vert \cdot \Vert$ and $\Vert \cdot \Vert'$ be two  adelic metrics on $L$. For every
place $v$ of $K$ there is a constant $C_v \geq 1$ such that $C_v=1$ for almost all $v$
and
$$\Vert \cdot \Vert'_v \leq C_v \Vert \cdot \Vert_v \, .$$ }
\end{lm}

 \medskip

\noindent {\bf Proof.} For every place $v$ of $K$,  the quotient $\Vert \cdot \Vert'_v/\Vert \cdot \Vert_v$ is 
a continuous function on   $X^{\rm an}_v$, and $X^{\rm an}_v$ is compact,
so there exists $C_v$ such that $\Vert \cdot \Vert'_v \leq C_v \Vert \cdot \Vert_v $.

The fact that  one can take $C_v=1$ for almost all $v$ follows from the fact that adelic metrics are coherent.
If $\mathcal O_K$ is the ring of integers in $K$, there exists a non empty open subset $U \subset {\rm Spec}(\mathcal O_K)$ and a model $(X_U,L_U)$ of $(X,L)$  over $U$
inducing both $\Vert \cdot \Vert_v$ and $\Vert \cdot \Vert'_v$ for every $v\in U$. 
Therefore one can assume that $C_v=1$ when $v\in U$.

\bigskip

\subsection{\\}  An {\it adelic vector space \it} over $K$ is a pair $\overline V = (V, \Vert \cdot \Vert)$
of a finite dimensional $K$-vector space $V$ and a collection $\Vert \cdot \Vert = \{\Vert \cdot \Vert_v\}$
of norms $\Vert \cdot \Vert_v$  on $V\otimes K_v$ for every place $v$ of $K$.
One assumes that $\Vert \cdot \Vert_v$ is non-archimedean when $v$ is finite and that there exists a non-empty
open subset $U \subset {\rm Spec}({\mathcal O}_K)$ and a projective module $V_U$ over ${\rm Spec}({\mathcal O}(U))$
such that $V= V_U \otimes_{{\mathcal O}(U)}K$ and $\Vert \cdot \Vert_v$ is induced by $V_U$
 when $v\in U$.

Let ${\cal M}(K)$ be the set of places of $K$. If ${\mathcal O}_v$ denotes the integers in $K_v$, let
$s\in V$ and $s\neq 0$. For almost all $v\in {\cal M}(K)$ the module $V_U \otimes {\mathcal O}_v$ is free and   $\Vert s \Vert_v=1$, hence  the following sum
$$h(s)= \sum_{v \in {\cal M}(K)}\log(\Vert s \Vert_v)$$
is finite.

\subsubsection{\\}

Let $X$ be a smooth projective variety over $K$ and let $\overline L = (L, \Vert \cdot \Vert)$ be an adelic line bundle
on $X$. If $v$ is a place of $K$ one defines as follows a norm $\Vert \cdot \Vert_{\sup,v}$ on $H^0(X,L)\otimes_K K_v$.
If $s\in H^0(X,L)\otimes_K K_v$
$$ \Vert s \Vert_{\sup,v} = \sup_{x\in X(\overline K_v)} \Vert s \Vert_v(x) \, .$$
The pair $H^0(X,L)_{\sup}=(H^0(X,L), \{\Vert \cdot \Vert_{\sup,v}\})$ is adelic (\cite{Z2},   (1.6)) . More precisely:

\begin{prop}\label{G} (Gabber)

 {\it Assume that $U\subset {\rm Spec}(\mathcal O_K)$ is a non-empty
subset and that $f: X_U  \to U $ is a flat projective morphism with reduced fibers and generic fiber
$X \to {\rm Spec}(K)$. Let $L_U$ be a line bundle on $X_U$ and let $L= L_U \otimes K$. When $v\in U$ we equip
$L\otimes K_v$ with the norm $\Vert \cdot \Vert_v$ induced by $L_U$. Then the sup-norm
$\Vert \cdot \Vert_{\sup,v}$ on $H^0(X,L)\otimes_K K_v$ is induced by $H^0(X_U,L_U)$.

}
\end{prop}

\subsubsection{\\} {\bf Proof.} Let $\Vert \cdot \Vert_{m,v}$ be the norm on $H^0(X,L)\otimes_K K_v$ 
induced by $H^0(X_U,L_U)$ : if $s\in H^0(X,L)\otimes_K K_v$ ,
$$\Vert s \Vert_{m,v} = \inf_{a\in \overline K_v} \{ |a| / s = a s_0, s_0 \in H^0(X_U,L_U) \otimes {\overline  {\mathcal O}_v} \} \, ,$$
where ${\overline  {\mathcal O}_v }$ are the integers in ${\overline  K_v }$.
We  first show that
$\Vert s \Vert_{\sup,v} \leq \Vert s \Vert_{m,v}$. By the valuative criterion of properness
$X({\overline  K_v }) = X_U ({\overline  {\mathcal O}_v )}$. Therefore, if $x \in X({\overline  {K}_v })$ we have
$\Vert s_0(x) \Vert_v \leq 1$. This implies $\Vert s(x) \Vert_v \leq |a|$, hence
$\Vert s \Vert_{\sup,v} \leq  |a| $ and  $\Vert s \Vert_{\sup,v} \leq \Vert s \Vert_{m,v}$.

\subsubsection{\\}

To show that $\Vert s \Vert_{m,v} \leq \Vert s \Vert_{\sup,v}$, let $\pi$ be a uniformizer of $ {\mathcal O}_v$ and let
$\overline {(k(v))}$ be an algebraic closure of $k(v)=  {\mathcal O}_v/\pi$. There exists a power $\pi^n$ of $\pi$
such that ${\widetilde {s}} = \pi^n s$ lies in $H^0(X_U,L_U) \otimes {\overline  {\mathcal O}_v}  $ and
$\widetilde {s}$ does not vanish on the special fiber $X_0$ of $f$ at $v$. There is thus a closed point $x_0 \in X_0$ such that
$\widetilde {s} (x_0) \neq 0 $. The specialization map from $X(\overline {K_v})$ to $X_0(\overline {k(v)})$
is surjective, so $x_0$ lifts to a point $\widetilde {x} \in X_U ({\overline  {\mathcal O}_v })$.  Since
$\widetilde {s} (x_0) \neq 0 $, $\Vert \widetilde s (\widetilde x) \Vert_v =1$ and $ \Vert \widetilde s \Vert_{\sup,v} \geq 1$.

On the other hand, since ${\widetilde {s}} $ lies in $H^0(X_U,L_U) \otimes {\overline  {\mathcal O}_v}  $
we have $\Vert \widetilde s \Vert_{m,v} \leq 1$. So $\Vert \widetilde s \Vert_{m,v} \leq \Vert \widetilde  s \Vert_{\sup,v}$
and $\Vert s \Vert_{m,v} \leq \Vert s \Vert_{\sup,v}$. This ends the proof of Proposition~\ref{G}.

\subsection{ \\ }  Let  $Y$ be a smooth, irreducible  and projective curve over $K$. If $F$ is a finite extension 
of $K$ and $P \in Y(F)$, we let $u : {\rm Spec} (F) \to Y$ be the map defined by $P$.  

\bigskip

\begin{prop}\label{prop1} {\it Let $\widehat{x} \in \widehat{\rm Pic} (Y)_{\mathbb Q}$ and let $x \in {\rm Pic} (Y)_{\mathbb Q}$ be the image of $\widehat x$. Assume that $\deg_K (x) < 0$. Then
$$
\widehat{\deg}_F (u^* (\widehat x)) \leq 0(d)
$$
where $d = [F:{\mathbb Q}]$.}
\end{prop}

\medskip

\noindent {\bf Proof.} We can assume that $ - \widehat x \in \widehat{\rm Pic} (Y)$ is the class
 of an adelic line bundle $\overline L$ on $Y$. 

The line bundle $u^*(L)$ on ${\rm Spec}(F)$ is a one-dimensional vector space $\Lambda$
over $F$:
if $\eta$ is the point of ${\rm Spec}(F)$, $\Lambda = u^*(L)_\eta$.
When $w$ is a valuation of $F$, we let $F_w$ be the completion of
$F$ at $w$, and we let $\overline{F_w}$ be an algebraic closure of $F_w$. The metric
$u^*(\Vert \cdot \Vert)$ on $u^*(L)$ is a collection of $\overline{F_w}$-norms
$u^*(\Vert \cdot \Vert)_w$ on $\Lambda \otimes_F \overline{F_w}$. They are defined as follows.
 Let $v$ be the restriction 
of $w$ to $K$.
Any element in  $ \Lambda \otimes_F \overline{F_w}$ is of the form $u^*(m)\otimes \alpha$
where  $m\in L_{u( \eta)}$ 
and  $\alpha \in {\overline K_v} = {\overline F_w}$ (see  Lemma~\ref{lemma100}  i) below).
We define
$$u^*(\Vert u^*(m)\otimes \alpha \Vert)_w = 
 \vert \alpha \vert_w \Vert m \Vert_v \, . \leqno (1) $$

Let us see that $u^*(\Vert \cdot \Vert)$ is adelic over $F$ in the sense of § 1.3.2. First,
since ${\rm Gal}(\overline{F_w}/F_w)$ is contained in  ${\rm Gal}(\overline{K_v}/K_v)$
the norm $u^*(\Vert \cdot \Vert)_w$ is invariant under  ${\rm Gal}(\overline{F_w}/F_w)$.
Next we indicate why 
$u^*(\Vert \cdot \Vert)$  is coherent. 
The metric $\Vert \cdot \Vert$ on $L$
being coherent there exists a non-empty open subset
$U\subset {\rm Spec}(\mathcal O_K)$ and a model $(Y_U,L_U)$ of $(Y,L)$  over $U$
inducing $\Vert \cdot \Vert$.  Denote by $V$ the inverse image of $U$
in ${\rm Spec}(\mathcal O_{F})$ and let $A$ (resp. $B$) be the ring of functions on $U$
(resp. $V$). 
 By the valuative criterion of properness the point $u( \eta)\in Y( F)$
gives a morphism $\widetilde u: V\to Y_U$. We let  $L_V= \widetilde u^*(L_U)$ on $V$.
The Lemma~\ref{lemma100} ii) below shows that $B\otimes_A  K = F$
and we get 
$$L_V \otimes_B F = L_V\otimes_B B\otimes_A K=\widetilde u^*(L_U) \otimes_A K =
u^*(L_U \otimes _A K) = u^*(L) \, ,$$
so $(V,L_V)$ is a model of $({\rm Spec}(F), u^*(L))$. It follows from
the definition (1) of the pull-back norms  
 that, for every place $w\in V$,
this model induces the norm $u^*(\Vert \cdot \Vert)_w$.
Therefore the metric $u^*(\Vert \cdot \Vert)$
 is adelic over $F$.
We denote by $u^*(\overline L) = -u^*(\widehat x)$ the class of 
 $(u^*(L), u^*(\Vert \cdot \Vert))$
in $\widehat{\rm Pic} ({\rm Spec}(F))$.

Since $\deg_K (L) > 0$, $L$ is an ample line bundle on $Y$.
 A power of $L$ is very ample: there is a closed immersion $i : Y \to {\mathbb P}_K^N$ and an integer $n \geq 1$ such that $i^* ({\cal O} (1)) = n L$. If $\Vert \cdot \Vert_{{\cal O}(1)}$ is the Fubini-Study metric at archimedean
places and the standard norm at finite places, the $n$-th root of $i^* (\Vert \cdot \Vert_{{\cal O}(1)})$ is
 a non-negative adelic metric $\Vert \cdot \Vert'$ on $L$ and thus $\widehat{\deg}_F (u^*(L , \Vert \cdot \Vert'))\geq 0$.

 By Lemma~\ref{lm1}  we can compare the adelic metrics $\Vert \cdot \Vert$ and $\Vert \cdot \Vert'$ on $L$.
If $w$ is a place of $F$ restricting to $v$ on $K$,  if $m\in u^*(L)_{ \eta}$,
and if $\alpha \in {\overline F_w}$, $m\neq 0$, $\alpha \neq 0$, we deduce from
(1) that
 $$ |\log(u^*( \Vert m\otimes \alpha \Vert)_w) - \log(u^*( \Vert m\otimes \alpha \Vert ')_w)| \leq \log(C_v) \, ,$$
 where $C_v \geq 1$ and  $C_v=1$
 for almost all $v$. Therefore
$$
\widehat\deg_F (u^* (\widehat x)) = - \widehat\deg_F (u^*(L,\Vert \cdot \Vert))\leq d \sum_v \log(C_v)  \le  0(d) 
\, .
$$

\bigskip

\begin{lm}\label{lemma100} Let $F$ be a number field containing $K$.

i) Let $w$ be a finite place of $F$ and let $v$ be its restriction to $K$. Denote by $K_v$ and $F_w$ the completions and by $\overline K_v$ (resp. $\overline F_w$) an algebraic closure of $K_v$ (resp. $F_w$). There exists a field isomorphism $\overline K_v \to \overline F_w$ which preserves the norms.

ii) Let $A$ (resp. $B$) be a domain with field of fractions $K$ (resp. $F$).
If the inclusion $K\subset F$ sends $A$ into  $B$ we have $B\otimes_A K=F$.
\end{lm}

\noindent {\bf Proof.}    
i) Since $K_v \subset F_w \subset \overline F_w$, the field $\overline F_w$ is an algebraic closure of $K_v$, so it is isomorphic to $\overline K_v$. Since $K_v$ is complete the extension of the valuation $v$
to the fields $\overline K_v$ and $\overline F_w$ is unique (and given by the norm).

ii) The ring $B\otimes_A K$ is obtained from $B$
by inverting the nonzero elements of $ A$, so it can be viewed as a subring of the field of fractions $F$ of $B$. Furthermore the map $B\otimes_A K \to F$ is surjective: let $a/b \in F$
with $a\in B$ and 
$b\in B-{\{0 \}}$. The element $b\in F$ satisfies a non-trivial equation
$a_0 b^n+a_1 b^{n-1}+ \cdots +a_n=0$ with $a_i \in A$, from which we deduce 
that there exists $c\in B-{\{0\}}$ with $bc \in A $. Therefore
$a/b = (ac)/(bc)$ lies in $B\otimes_A K$.

\subsection{\\}
We keep the notation of the previous paragraph. Let ${\overline M}$ be an adelic line bundle on $Y$ with $M$
an ample line bundle. 

\begin{cor}\label{cor1} {\it Let $\widehat{x} \in \widehat{\rm Pic} (Y)_{\mathbb Q}$ and let $x \in {\rm Pic} (Y)_{\mathbb Q}$ be the image of $\widehat x$. Assume that $\deg_K (x) = 0$. Then, for every positive rational
number $\varepsilon$, one has
$$
\vert\widehat{\deg}_F (u^* (\widehat x))\vert \leq \varepsilon \widehat{\deg}_F (u^* (\overline M)) + 0(d) \, .
$$}
\end{cor}

\medskip

\noindent {\bf Proof.} Apply Proposition~\ref{prop1} to $\widehat x - \varepsilon \overline M$
and $-\widehat x -\varepsilon \overline M$. Note that the implied constant in $0(d)$ depends on $\varepsilon$.

 \bigskip

\subsection {\\}

Let $\overline M$ be an adelic line bundle on ${\rm  Spec}(K)$.
Assume that $M$ has dimension $2$ and that $\overline M$ has a model
on the whole of ${\rm Spec}({\mathcal O}_K)$. 
 Minkowski's theorem on successive minima implies that   there exist two sections $s_1$ and $s_2$
in $M$ which are linearly independent over $K$ and such that 
$$ h( s_1)+ h( s_2) \leq   -  \widehat \deg_K(\overline M) +
 \delta(K)  +  d_K\, A \,  ,  \leqno(2)$$
where $A$ is an absolute  constant  (see \cite{BV} Th.~3  and \cite{BO} Th.~2).

\bigskip

\section{An inequality}

\bigskip

\subsection{\\}

\subsubsection {\\}
Let $X$ be a smooth projective curve of positive genus $g$ over $K$ and let $v$ be a place of $K$.

When $v$ is archimedean and given by a complex embedding $\sigma$ of $K$ we let $\mu_v$ be the 
following probability measure on $X(\overline K_v) = X_\sigma$ 
(\cite{A}, §4). If $\alpha$ and $\beta$ are holomorphic one-forms on $X_\sigma$
we let
$$(\alpha,\beta)= \frac {i}{2} \int_{X_\sigma} \alpha \wedge {\overline \beta} \, .$$
For any orthonormal basis $\omega_1, \cdots, \omega_g$ of $\Omega^1(X_\sigma]$
we define
$$\mu_v = \frac{i}{2g} \sum_{j=1}^g\omega_j\wedge .\overline \omega_j  \, .$$

\medskip

When $v$ is finite we let $\mu_v$ be the probability measure on $X_v^{\rm an}$
defined in \cite {Z1}, Lemma 3.7, and  \cite {Y} , Prop. A.5.

\subsubsection{\\} 
Let $Y$ be a smooth projective variety over $K$ and let $f:X \to Y$ be a smooth connected curve
of positive genus over $Y$. If $\overline L$ is an adelic line bundle over $X$ we say that
$\overline L$ is $admissible$ when $\overline L$ is integrable and, for any place $v$ of $K$
and any point $y\in Y({\overline K_v})$, the curvature of the restriction of $\overline L$
to the geometric fiber of $f$ at $y$ is proportional to $\mu_v$.
It is shown in \cite{Y}, Prop. 2.6, that any line bundle on $X$ has an admissible metric

\medskip

\subsection{\\}

  Let  $X$ be a smooth, geometrically irreducible and projective curve of positive genus over $K$. Assume $\overline L_1$ and $\overline L_2$
are two ample line bundles  over $X$ equipped with  admissible     adelic metrics.
If $i=1$ or $2$, let $d_i$ be the degree of $L_i$
and let
 $s_i \in H^0(X,L_i)$ be a non-trivial global section of $L_i$ on $X$.
 We assume that the divisors of $s_1$ and $s_2$
do not meet.
Let 
$$h(s_i)= \sum_ {v\in {\cal M}(K)}\log(\Vert s_i \Vert_{{\rm sup},v}) \, .$$

\medskip

\begin{prop}\label{prop12}

\bigskip
The following inequality holds:
$$
- {\overline L_1}\cdot{\overline L_2}  \leq d_2 h(s_1) + d_1 h(s_2)   \, .
$$

\end{prop}.

\noindent {\bf Proof.}
 Let $<s_1,s_2> \neq 0$ be the element of  $<L_{1},L_{2}>$ associated to
 $s_1$ and $s_2$.
For any place $v$ of $K$ we equip the Deligne pairing $<L_{1,v},L_{2,v}>$ of $L_{1,v}$ and $L_{2,v}$
with the metric
defined in \cite{De}(6.3.2) and \cite{YZ} (4.1.1) when $v$ is archimedean, and 
\cite{Z1} p.177 when $v$ is finite.
If $g_i = - \log \Vert s_i \Vert_v$ we have
$$- \log (\vert <s_1,s_2> \vert_v)  = \int_{X_v^{\rm an }} (g_1 {\delta}_{Z_2} + c_1(\overline{L_1})g_2) \, , \leqno [3)$$
where ${\delta}_{Z_2}$ is the current of integration on $Z_2 = {\rm div} (s_2)$, and $c_1(\overline{L_1})$
is the first Chern form of $(L_{1,v}, \Vert \cdot \Vert_v)$
(see \cite  {Z1} p.177 l.18,  (2.7) and (4.5)). 

The divisor $Z_2 = \sum_\alpha n_\alpha P_\alpha$ is effective therefore
$$\int_{X_v^{\rm an }} g_1 {\delta}_{Z_2} = \sum_\alpha - n_\alpha  \log \Vert s_1 \Vert_v (P_\alpha) \geq - d_2 \log \Vert s_1 \Vert_{\sup, v } \, . \leqno(4)$$
Since the metric on $\overline L_1$ is admissible, the first Chern form
$c_1(\overline{L_1}) = d_1 \mu_v$ is non-negative and
$$\int_{X_v^{\rm an }}c_1(\overline{L_1})g_2 \geq -d_1 \log \Vert s_2 \Vert_{\sup, v } \, . \leqno [5)$$
From (3), (4) and (5) we deduce that
$$- \log (\vert <s_1,s_2> \vert_v) \geq - d_2 \log \Vert s_1 \Vert_{\sup, v }  -d_1 \log \Vert s_2 \Vert_{\sup, v }  \, . \leqno [6)$$

By definition, the arithmetic intersection number $\overline L_1 \cdot \overline L_2$ is the sum 
\noindent $-\sum_v \log( \vert\langle s_1 ,s_2 \rangle \vert_v)$. Therefore, by (6),  we get
$$
- {\overline L_1}\cdot{\overline L_2}  \leq d_2 h(s_1) + d_1 h(s_2)   \, .
$$

\bigskip

\section{On some curves of genus three and six}

\subsection{ \ } Let $\overline K$ be an algebraically closure  of  $K$,  and let $Y_0$ be a smooth projective connected curve of genus $g_0 = 2$ on $\overline K$. Let $J_0$ be the jacobian of $Y_0$ and
\be
\nonumber
\chi : Y_1 \to Y_0
\ee
a non trivial \'etale map of degree two, classified by a class
\be
\nonumber
\eta \in H^1_{\rm {et}} (Y_0, {\mathbb Z} / 2) = J_0 [2] \, .
\ee
Let $P_0$ be a point in $Y_0(\overline K)$ and $\{ P_1 , P_2 \} = \chi^{-1} (P_0)$. We assume     given a  
proper surjective morphism of degree two
\be
\nonumber
\theta : V \to Y_1
\ee
which is ramified above $P_1$ and $P_2$. 
By
 \cite {BPV}, I 17,
 the choice of $V$ corresponds to the choice of a square root $\lambda$ of $ {\mathcal O}_{Y_1}(P_1 + P_2)$:
 $$ 2\lambda = P_1+P_2 \, . $$
 Note that, by the Riemann-Hurwitz formula, $Y_1$ has genus $g_{Y_1}= 3$ and $V$ has genus $g_V=6$.
 Let $\omega_1$ be the sheaf of differentials on $Y_1$.

\bigskip

\noindent 
\begin{prop}\label{prop3}
\item[i)] The curve $Y_1$ is hyperelliptic.
\item[ii)] Assume $W$ is a Weierstrass point on $Y_1$. Then
$4 {\mathcal O}_{Y_1}(W)=\omega_1$.
\item[iii)] Assume  $P_1$ is a Weierstrass point on $Y_1$. Then $h^0(P_1+\lambda)=0 $.
\end{prop}

\bigskip

\subsection  { \ }
 {\bf Proof.}   Let us prove Proposition~\ref{prop3}  i). 
Since $\{ P_1 , P_2 \} = \chi^{-1} (P_0)$ we have
$$
h^0 (Y_1 , P_1 + P_2) = h^0 (Y_1 , \chi^* (P_0)) = h^0 (Y_0 , \chi_* \, \chi^* (P_0)) \, .
$$
Since
$$
\chi_* ({\mathcal O}_{Y_1}) ={\mathcal O}_{Y_0} \oplus  \eta
$$
(\cite{BPV}, I  (17.2)) 
we get
$$
h^0 (Y_1 , P_1 + P_2) = h^0 (Y_0 , P_0) + h^0 (Y_0 , P_0+\eta) \, .
$$

The divisor $\eta$ is the difference of two points of $Y_0$.
To prove that, we consider $Y_0$ as an imaginary hyperelliptic curve,
and we let $\infty$ be the point at infinity.
Let us write $\eta$ in reduced form (\cite{MYZ}, Theorem 6.1 \footnote{On an imaginary hyperelliptic curve of genus
$g$ a divisor is called {\it reduced\it} if it is of the form $\sum_{i=1}^{i=k} P_i- k\infty$
with $k\leq g$, $P_i\neq \infty$ if $i =1 \cdots k$, and $P_i + P_j \neq 2\infty$ when $i\neq j$.}):
$$\eta = A_1+A_2-2\infty \,  ,$$
 or $\eta = A_1 - \infty$.
If $P_0 = 2\infty - A_2$ is the conjugate 
 of $A_2$ on the hyperelliptic curve $Y_0$ or if $P_0= \infty$, we get
$$\eta = A_1-P_0\, .$$

 Therefore 
$$
h^0 (Y_1 , P_1 + P_2) = h^0 (Y_0 , P_0) + h^0 (Y_0 , A_1) = 2
$$
and $P_1 + P_2$ is a $g_2^1$ on $Y_1$. This proves i).

 If $W$ is a Weierstrass point, $2W$ is the hyperelliptic
divisor, hence $4 {\mathcal O}_{Y_1}(W)=\omega_1$. This proves ii).

\subsection{ \ }

\medskip

Let us prove Proposition~\ref{prop3} iii).
 We proceed by contradiction. If $h^0 (P_1 + \lambda)$ is positive there exist two points $P$ and $Q$ on $Y_1$ such that 
$$
P_1 + \lambda = P+Q \, . \leqno (7)
$$
Since $Y_1$ is hyperelliptic (see i)), we may think of it as an imaginary hyperelliptic
 curve with one point $\infty$ at infinity. Since $P_1$ is a Weierstrass point we can assume that $P_1=\infty $.
If $P$ is a point on $Y_1$, we let
$$
 \widetilde P = 2 \infty - P
$$
be its conjugate.

Multiplying  (7) by $2$ we deduce that
$$
 P_2 +  \widetilde P +  \widetilde Q  - 3 \infty=   P + Q - 2 \infty \, . \leqno (8)
$$
If $P\neq \infty$, $Q \neq \infty$ $P_2 \ne P$,  $P_2\ne Q$ and $P\ne \widetilde Q$ both sides  of (8) are reduced.
In that case, by the uniqueness of the reduced form of a divisor class of degree zero (\cite{MYZ}, Theorem 6.1), (8) 
cannot happen since the length of the right-hand side differs from the length of the left-hand side.

\smallskip

Assume that $P = \infty$. Then (8) becomes
$$
 P_2 +  \widetilde Q  - 2 \infty=    Q - \infty \, . \leqno (9)
$$
If $Q\ne \infty$ and $P_2 \ne Q$ both sides of (9) are reduced; this cannot happen since both sides do not have the same length.
If $Q=\infty$ we get $P_2=\infty $ which is impossible since $P_2 \ne P_1$.
If $P_2=Q$ we get $Q= \infty$, which is impossible for the same reason.

By permuting $P$ and $Q$ we see that it is also not possible that $Q=\infty$.

 Assume now that $P_2 = P$. From (8) we get
$$
 \widetilde Q - \infty = P_2 +Q - 2 \infty \, . \leqno (10)
$$
Again it cannot happen that both sides of (10) are reduced.
Therefore $Q=\infty$, $P_2=\infty$ or $P_2= \widetilde Q$.
In these three cases $P_2= \infty$, which is impossible.
Similarly it cannot happen that $P_2 = Q$.

Finally, if $P= \widetilde Q$ we get from (8) the equality
$P_2=\infty$, which is impossible.
This ends the proof of Proposition~\ref{prop3}.

\bigskip

\subsection {\\}

\medskip

\begin{prop} (Voisin) \label{prop13}
\begin{enumerate}
\item[i)] We have 
$h^0(P_1+P_2)=2$ if and only if  there exists a point $P'_0$ on $Y_0$ such that
$
\eta = P'_0 - P_0 
$. 
\item[ii)] We have $h^0 (Y_1 , \lambda) = 0$, except when $\lambda$ is 
a Weierstrass point on $Y_1$, in which case
$$h^0 (Y_1 , \lambda) = 1  \,  $$
and  $h^0(P_1+P_2)=2$.
\item[iii)] If $h^0(P_1+P_2)=2 $ we say that  $P_1$ is an {\em exceptional point}. There are 4 exceptional points on  $Y_1$. The point $P_1\in Y_1$ is exceptional if and only if $P_1$ and $P_2$ are conjugate by the hyperelliptic involution of $Y_1$.
\end{enumerate} 
\end{prop}

\bigskip

\subsubsection  { \ }

 {\bf Proof.}   Let us prove Proposition~\ref{prop13}  i). 
We showed  in the proof of Proposition~\ref{prop3}  i) that 
$$
h^0 (Y_1 , P_1 + P_2) = h^0 (Y_0 , P_0) + h^0 (Y_0 , P_0+\eta) \, .
$$
This proves the  assertion.

\bigskip

\subsubsection{ \ }

To prove Proposition~\ref{prop13} ii) we note that, since $\lambda$ has degree one, $h^0 (\lambda) = 0$ except when there exists a point $Q$ on $Y_1$ such that $\lambda = Q$, hence
$$
P_1 + P_2 = 2Q \, . \leqno (11)
$$
When (11) holds, the linear system $\vert 2\lambda \vert$ contains $P_1 + P_2$ and $2Q$. Since $P_1 \ne P_2$, this implies that $2\lambda$ is the hyperelliptic divisor on $Y_1$. Conversely, if $P_1 + P_2$ is the hyperelliptic divisor on $Y_1$, $P_1 + P_2 = 2W$ for every Weierstrass point $W$ on $Y_1$ and if (11) holds $Q=2W-Q$ is a Weierstrass  point.  It follows that $h^0(P_1+P_2)=h^0(2W)=2$ .
 This proves ii).

\bigskip 

\subsubsection{ \ }

Next we prove Proposition~\ref{prop13} iii).
If $P$ and $P'$ are two Weierstrass points on $Y_0$, their difference $P-P'$ has order two in the Jacobian $J_0$ of $Y_0$.
If $P$, $P'$, $R$ and $R'$ are points on $Y_0$ such that $P \neq P'$
and 
$$P-P'=R-R' \, ,$$
we get
$$P+R'-2\infty = P'+R-2\infty \, ,$$
hence, by the uniqueness of the reduced form of a divisor class of degree zero
(\cite{MYZ}, Theorem 6.1), $P=R$ and $P'=R'$ or $P=2\infty-R'$ and $P'=2\infty-R$.
Therefore the map sending a pair of distinct Weierstrass points to their difference is injective.
 There exist $2g_0 + 2 = 6$ Weierstrass points on $Y_0$, and the 2-torsion subgroup $J_0[2]$ of $J_0$ has order $2^4 = 16$. By counting we see that, if $\eta' \ne 0$ lies in $J_0[2]$, there exists a unique pair $\{P,P'\}$ of Weierstrass points on $Y_0$ such that
$$
\eta' = P-P' = P'-P \, .
$$

In particular, there is a unique pair of Weierstrass points $\{P_0,P'_0\}$
on $Y_0$ such that 
$$\eta = P_0-P'_0 = P'_0 - P_0 \, .$$ 
Furthermore, if $\eta$ is the difference of two points, these are Weierstrass points.
Therefore, by i) and ii), $P_1$ is exceptional if and only if 
$\chi(P_1)=P_0$ or $P'_0$.
Since $\chi$ has degree two,  there are four exceptional points on $Y_1$.

Finally, by (ii), the point $P_1$ is exceptional if and only if, for any Weierstrass point $W$, $P_2=2W-P_1$, i.e. $P_1$ and $P_2$ are conjugate
by the hyperelliptic involution of $Y_1$. 

\bigskip
\subsection{\\}

We keep the notation of Proposition~\ref{prop13}.

\begin{prop}\label{prop15}
Assume that $P_1$ is a  Weierstrass point . There exists  a divisor $D$ of degree $1$ on $Y_1$
such that $2D=P_1+P_2$ and,  if $\rho: V\to Y_1$ is the covering of $Y_1$ defined by $D$,
 there exists no point $Q$
on $V$ such that
$$h^0(V,\rho^*(P_1+P_2)+Q) = 2 \, . \leqno (12) $$
\end{prop}

\bigskip

 {\bf Proof.} (Voisin)  

Let $\Theta$ be the Theta divisor of $Y_1$, i.e. the locus of classes of effective divisors
of degree two on $Y_1$. The Proposition~\ref{prop15} follows from
the following lemma.

\begin{lm}\label{lm15}
Let $P_1\in Y_1(\overline K)$ and let $P_2$ be the image of $P_1$
by the involution of $Y_1$ over $Y_0$. There exist two points $R_1$ and $R_2$ in $Y_1(\overline K)$
such that, if $D$ is a divisor of degree $1$ on $Y_1$
with $2D=P_1+P_2$, if $h^0(D+P_1)=h^0(D+P_2)=0$,  and if there exists   a point $Q$ on $V$
such that (12) holds, then $D+R_1\in \Theta$ or $D+R_2\in \Theta$.
\end{lm}

\subsubsection{\ }

Let us prove that Lemma~\ref{lm15} implies Proposition~\ref{prop15}. Let $J_1$ be the jacobian of
$Y_1$. Fix a divisor $D_0$ of degree $1$ such that $2D_0 = P_1+P_2$.
Any divisor $D$ of degree $1$ satisfying $2D=P_1+P_2$
is of the form $D=D_0+\eta$ with $\eta \in J_1(\overline K)$ and $2\eta=0$.
Since $P_1$ is a Weierstrass point,  Proposition~\ref{prop3} iii) implies that  $h^0(D+P_1) = h^0(D+P_2)=0$.
A translate of $\Theta$ in $J_1$ contains at most $29$ points of order $2$
(\cite{MU}, p. 3.105).
There are $64$ points of order $2$ in $J_1$. Since $64>2\times 29$
we may choose $\eta$ such that $\eta + D_0 +R_i \notin \Theta$ for $i=1$ and $2$.
By  Lemma~\ref{lm15} there does  not exist any point $Q$ on $V$ satisfying the equation (12).

\subsubsection{\ }

 Let $\lambda = {\mathcal O}_{Y_1} (D)$. 
Since $V$ has genus $6$, Serrre duality and the Riemann-Roch theorem
imply
$$h^0(V,\omega_V-\rho^*(P_1+P_2)-Q) = h^1(V,\rho^*(P_1+P_2)+Q) $$
$$= h^0(V,\rho^*(P_1+P_2)+Q)-5+6-1= h^0(V,\rho^*(P_1+P_2)+Q) \, .$$
Let $Q_1$ (resp. $Q_2$) be the inverse image of $P_1$ (resp. $P_2$) in $V$.
By the Riemann-Hurwitz theorem $\omega_V= \rho^*(\omega_1) +Q_1+Q_2$. 
Therefore (12) is equivalent to the equality
$$ h^0(V,\rho^*(\omega_1-\lambda)-Q)=2.$$

\begin{lm} \label{lm16} 
$$h^0(V,\rho^*(\omega_1-\lambda)) = 2 \, .$$
\end{lm}

\noindent{\bf Proof.} Since $\rho_*({\mathcal O}_V) = {\mathcal O}_{Y_1} \oplus \lambda$
(\cite {BPV}, I, (17.2)) we get
$$h^0(V,\rho^*(\omega_1-\lambda)) = h^0(Y_1,\omega_1-\lambda) + h^0(Y_1,\omega_1-P_1-P_2) \, .$$
It is thus enough to show that $ h^0(Y_1,\omega_1-\lambda)=1$ and $h^0(Y_1,\omega_1-P_1-P_2) =1$.
Now $h^0(Y_1,\omega_1-\lambda)=h^0(Y_1,\lambda)-1+3-1=h^0(Y_1,\lambda)+1$.
Since $P_1$ is not exceptional $h^0(Y_1,\lambda)=0 $ (Proposition~\ref{prop3}, iii)), therefore
$h^0(Y_1,\omega_1-\lambda)=1$. Similarly
$h^0(Y_1,\omega_1-P_1-P_2) = h^0(P_1+P_2)-2+3-1=h^0(P_1+P_2)=1$
by
Proposition~\ref{prop3}, iii).

\subsubsection{\\}
Let us prove  Lemma~\ref{lm15}. By Lemma~\ref{lm16}, (12) is equivalent to the fact that
$Q$ is a base point of the linear system $|\rho^*(\omega_1-\lambda)|$ on $V$. 

Notice that $Q$ is not equal to $Q_1$ or $Q_2$. Indeed, the sections of $\omega_1-\lambda$ vanish
on $\rho(Q)$. But 
$$h^0(\omega_1-\lambda-P_1)=h^0(\lambda+P_1)=0$$
by assumption. Therefore $\rho(Q)\neq P_1$ and similarly $\rho(Q) \neq P_2$.

Let $\sigma : V \to V$ be the involution of $V$ over $Y_1$. Since $Q_1$ and $Q_2$ are the fixed points of
$\sigma$, the point $Q'=\sigma(Q)$ is another base point of $|\rho^*(\omega_1-\lambda)|$ on $V$. 
Let $R\in Y_1$ be the point such that $\rho^*(R)= Q+Q'$. By the argument in Lemma~\ref{lm16},  since
$h^0(\rho^*(\omega_1-\lambda-R))=2$, 
the sections in $ H^0(Y_1,\omega_1-\lambda)$ and those in $H^0(Y_1,\omega_1-P_1-P_2) $ vanish on $R$.
There are two points $R_1$ and $R_2$ on which the sections in $H^0(Y_1,\omega_1-P_1-P_2) $ vanish. The point  $R$ is equal to $R_1$ or $R_2$ and  $\omega_1-\lambda-R$ is effective, so
$$h^0(D+R)=h^0(\omega_1-\lambda-R) > 0 \, ,$$ i.e. $D+R \in \Theta$. This proves  Lemma~\ref{lm15}.

\section{Stable pointed curves}\label{sec6}

\bigskip

\subsection{ \ }
Fix an integer $n \geq 0$. Recall the following definition of Deligne and Mumford \cite{DM}. Let $Y$ be a scheme and $\phi : X \to Y$ a family of nodal curves over $Y$. Assume given sections $s_i : Y \to X$, $i = 1,\ldots , n$, with images $B_i = s_i (Y)$. These data define a {\it stable $n$-pointed curve} of genus $g$ over $Y$ when
\begin{enumerate}
\item[1)] The schemes $B_i$ are mutually disjoint;
\item[2)] The schemes $B_i$ are disjoint from the non-smooth locus of $\phi$;
\item[3)] All the geometric fibers of $\phi$ have genus $g$;
\item[4)] The sheaf $\omega_{X/Y} \left( \underset{i}{\sum} \, B_i \right)$ is $\phi$-ample.
\end{enumerate}

\medskip

If $Y$ is regular, it follows from 2) that each $B_i$ defines a line bundle ${\mathcal O}_X (B_i)$ equipped with a canonical section.
Indeed, since $B_i$ is contained in the regular locus, its ideal sheaf  is locally principal, generated by a non-zero divisor, and hence it
 determines a Cartier divisor.

Assume given a stable $n$-pointed curve $X$ over $Y$. Then, for every morphism $Y' \to Y$, the fiber product
$X'= X\times_Y Y'$ is a stable $n$-pointed curve over $Y'$.

Let $S$ be the spectrum of a Dedekind ring and $f : X \to S$, $f' : X' \to S$ two stable $n$-pointed curves over $S$. Assume that the generic fibers of $f$ and $f'$ are isomorphic (as pointed curves). Then there exists a unique isomorphism of $f$ and $f'$ over $S$ extending the generic one (\cite{AO}, Lemma~3.14).

\bigskip

\subsection{ \ }

We shall need the following extension result. Let $Y$ be a smooth projective connected curve over $K$ and $\overline Y$ a projective regular model of $Y$ over $S= {\rm Spec}({\mathcal O}_K)$.

\bigskip

 \begin{prop} \label{prop19}
Let $\phi: X\to Y$ be a stable smooth $2$-pointed curve of genus $g$ over $Y$. Then there exists an alteration
$a_S: Y_S \to {\overline Y}$ (i.e. $a_S$ is proper, surjective and generically finite)   inducing a finite \'etale map
$a: Y' \to Y$, and a stable $2$-pointed curve $\phi_S:  X_S \to Y_S$ of genus $g$ extending the curve 

\noindent $a_S^{-1}(\phi) : a_S^{-1}(X) \to a_S^{-1}(Y)=Y'$ .
\end{prop}

\medskip

\subsubsection{ \ } {\bf Proof.} (Abramovich)

Let $M_{g,2}$ (resp. ${\overline M}_{g,2}$) be the moduli space of smooth (resp. stable) $2$-pointed curves of genus $g$.
If $n>2$ is an integer we denote by ${\overline M}_{g,2}^n$ the moduli scheme
of stable $2$-pointed curves of genus $g$ equipped with a full level-$n$ structure. Assume that  $n=pq$
is the product of two primes bigger than $2$. According to \cite{Se1} and \cite{DM}, §5, 
${\overline M}_{g,2}^p$ (resp. ${\overline M}_{g,2}^q$) carries a family of stable $2$-pointed curves of genus $g$
over ${\rm Spec}({\mathbb Z}[1/p])$ (resp. ${\rm Spec}({\mathbb Z}[1/q])$). These families coincide over 
${\rm Spec}({\mathbb Q})$ so they define a familly over ${\rm Spec}({\mathbb Z})$

Since $Y$ carries a  smooth family of curves, there is a finite \'etale map $Y^n \to Y$ where these curves carry a full level-$n$ 
structure. This family defines a morphism $f: Y^n \to {\overline M}_{g,2}^n$. Let $Y_S$ be the Zariski closure
of the graph of $f$ in

\noindent  ${\overline Y} \times_{{\rm Spec}({\mathbb Z})} {\overline M}_{g,2}^n$.
Since $Y_S$ maps to $ {\overline M}_{g,2}^n$, it  carries a family of stable $2$-pointed curves of genus $g$.
The projection $a_S: Y_S \to {\overline Y}$ is an alteration which is finite and \'etale over $Y$.

\medskip

\subsubsection{ \ }

We shall denote $\phi' : X' \to Y'$ the map $a_S^{-1}(\phi)$.
We have two cartesian squares:

$$
\xymatrix{
X' \ar[d]_{\phi'} \ar[r]^v &X_S \ar[d]^{\phi_S} \\
Y' \ar[r] &Y_S \, .
}
$$
and

$$
\xymatrix{
X' \ar[d]_{\phi'} \ar[r] &X \ar[d]^{\phi} \\
Y' \ar[r]^a &Y \, .
}
$$

\bigskip

\section{The Kodaira fibration}

\subsection{ \ } Let $K_0$ be a number field and ${\mathcal C}$ a smooth, projective and geometrically 
connected curve over $K_0$, of genus $g_0=2$. As in \cite{MD} Proposition~1 (see also \cite{K} and \cite{BPV} p.~220),
 we shall define a finite extension $K$ of $K_0$, a surjective \'etale morphism
\be
\nonumber
 Y \to {\mathcal C} \underset{K_0}{\otimes} K
\ee
and a non-isotrivial,  smooth and proper map
\be
\nonumber
\phi  : X \to Y
\ee
with geometrically irreducible, one dimensional fibers of genus   $4g_0-2 = 6$.

\medskip

They are obtained as follows. Let $Y_0 = {\mathcal C} \underset{K_0}{\otimes} K$ and let 
$\chi : Y_1 \to Y_0$ be a non-trivial surjective \'etale map of degree two, which exists if $K$ is large enough.
Consider the divisor $\Gamma = Y_1 \underset{Y_0}{\times} Y_1$ in $Y_1 \times Y_1$. This divisor has two components:
$$
\Gamma = \Delta_1 \amalg \Delta_2 \, , 
$$
where $\Delta_1$ is the diagonal and $\Delta_2$ is the  graph of the involution $\nu$ of $Y_1$ over $Y_0$.
We also denote by $\Delta_0$ the diagonal in $Y_0\times Y_0$.
Let $J_1$ be the jacobian of $Y_1$ and
let ${\mathcal P}$ be a Poincar\'e line bundle on $Y_1 \times J_1$. There is a morphism
$$
\psi_1 : Y_1 \longrightarrow J_1
$$
and a line bundle $L$ on $Y_1$ such that 
$$
{\mathcal O}_{Y_1 \times Y_1} (\Gamma - 2 \Delta_1) = (1 \times \psi_1)^* ({\mathcal P}) +\pi_ {2}^* (L) \, , 
$$
where $\pi_2 : Y_1 \times Y_1 \to Y_1$ is the second projection.
 We let $Y$ be a connected component of the fiber product $Z$ of $\psi_1$ with the multiplication by $2$ in $J_1$. So  we get a  commutative diagram
$$
\xymatrix{
Y\ar[r]^-{\psi_2} \ar[d]_-\psi &J_1 \ar[d]^-{\times 2} \\
Y_1 \ar[r]^-{\psi_1} &J_1 \, ,
} 
$$
with $\psi$ finite and \'etale.

\smallskip

Let $\widetilde\Gamma$ be the inverse image of $\Gamma$ in $Y_1\times Y$. If $K$ is large enough one can show \cite{MD} that there exists a line bundle ${\mathcal L}$ on $Y_1\times Y$ such that
$$
2 {\mathcal L} = {\mathcal O}_{Y_1\times Y} (\widetilde \Gamma) \, .
$$
A standard construction (see e.g. \cite{BPV} I, 17) associates to ${\mathcal L}$ a degree two finite surjective map
$$
\rho : X \to Y_1\times Y \, ,
$$
branched along  $\widetilde\Gamma$ and
such that
$$
\rho_{*} ({\mathcal O}_{X}) = {\mathcal O}_{Y_1\times Y} \oplus {\mathcal L}^* \, , \leqno   (13)
$$
where ${\mathcal L}^*= - {\mathcal L}$ is the dual of ${\mathcal L}$. We let
$$
\phi: X \to Y
$$
be the composite of $\rho$ with the projection $p$ of $Y_1\times Y$ onto the second factor. The map $\phi$
is smooth.
Also, when $i=1$ or $2$ (resp. $i=0$), we let $ D_i$ ( resp. $D_0$)  be the inverse image of $\Delta_i$ (resp. $\Delta_0$) in $Y_1\times Y$ (resp. in $Y_0\times Y$), and
$E_i= \rho^{-1}( D_i)\subset X$ if $i=1$ or $2$.  We then call $A_i={\mathcal O}_{Y_1\times Y} ( D_i)$.
and $B_i= {\mathcal O}_{X} (E_i)$ the corresponding line bundles.
Notice that both maps $E_i \to Y$  are isomorphisms
\medskip

\subsection{ \ }

Let $ \overline K$ be  an algebraic closure of $K$,  and let
  ${\overline y} \in Y(\overline K )$  be a geometric point. We  let $V=X_{\overline y}$  be the  geometric fiber at ${\overline y}$ of the projection $\phi : X \to Y$.
 The smooth  curve $V$ is a degree two finite surjective cover of the geometric fiber $Y_1\otimes_K{\overline K}$
  at ${\overline y}$ of the second projection $Y_1\times Y \to Y$. The proper surjective map
$$
\theta : V \longrightarrow Y_1\otimes _K\overline K$$
has branch-locus the points $P_1$ and $P_2$, where
$ P_i $
is the fiber at ${\overline y}$ of the map
$D_i \to Y$,
$i = 1,2$.  The projections of $P_1$ and $P_2$ to $  (Y_0 \times Y)_{\overline y}$
are equal to the point $P_0$, fiber at ${\overline y}$ of the map $D_0 \to Y$.

If  $
\lambda = {\mathcal L}_{\overline y}
$ is the restriction of ${\mathcal L}$ to $Y_1\otimes_K \overline K$, the ramified proper surjective map $ \theta $  is  associated to  the square root
$\lambda $
of $P_1 + P_2$.  We are thus in the situation of \S4.1
(see \cite{MD} loc.cit.). We let $\{Q_i\} = \theta^{-1} (P_i) $  if $ i = 1,2$. 

Note that
$$  P_1 = \psi({\overline y}) \, $$
and
$$ P_2 = (\nu \circ \psi)({\overline y})\, ,$$
where $\nu$ is the involution of $Y_1$ over $Y_0$.

\bigskip

\subsection{\\}
We shall now compute some (classical) intersection numbers on the surfaces $X$ and $Y_1\times Y$.
  As in \S 1.2, if $L$ and $M$ are line bundles on $X$ or  $Y_1\times Y$,
we consider  the Deligne pairing  $\langle L,M \rangle$ on $Y$ and $L \cdot M$ its degree.
We let $\phi_1: X \to Y_1$ be the composite of $\rho$ with the projection $p_1 : Y_1\times Y \to Y_1$ onto the first factor.
 Furthermore, we denote  by $\omega=\omega_{X/Y}$ the sheaf of relative differentials of $X$ over $Y$, and by $k={\rm{ deg}}(Y/Y_1) = \deg(\psi)$ the degree of $Y$ over $Y_1$.  Finally, we let $W$ be the line bundle associated
 to a Weierstrass point on $Y_1$. 
\bigskip

\begin{prop}\label{prop5}
{ \ }
One has
\be
\nonumber
\begin{matrix}
&A_1\cdot A_2 =0 \, , \hfill \\
&A_{1}^2 = A_2^2 =  -4 k  \, , \hfill \\
&B_{1}^2 = B_2^2 =  -2 k  \, , \hfill \\
&p_1^*(W) \cdot A_1 =p_1^*(W)\cdot A_2   = k \,  ,  \hfill\\
&\phi_1^*(W) \cdot B_1 =\phi_1^*(W)\cdot B_2   = k \,  ,  \hfill\\
&\phi_1^*(W)^2=0  \,  ,  \hfill \\
&\omega^2 = 12  k  \, . \hfill
\end{matrix}
\ee

 \end{prop}

\bigskip

\noindent{ \bf  Proof. } 
First, since  $E_1$ and $E_2$ do not intersect,
$A_1 \cdot A_2 = B_1 \cdot B_2=0$ .
Next, by the adjunction formula, we have
$$
\mathcal O_{Y_1\times Y_1}( \Delta_1)^2 = -2 (g_{Y_1}-1) \, ,
$$
where $g_{Y_1}=3$ is the genus of $Y_1$.  Since 
$$A_1=(1\times \psi)^*(\mathcal O_{Y_1\times Y_1}( \Delta_1))$$
 by the projection formula we get
 $$A_{1}^2 = ((1\times\psi)^*(\mathcal O_{Y_1\times Y_1}( \Delta_1)))^2 =
  k(\mathcal O_{Y_1\times Y_1}( \Delta_1))^2=-4k \, .$$ 
  Let $\nu$ be the involution of $Y_1$ over $Y_0$.
  Since $(\nu \times 1)^*(A_1)=A_2$ we have $A_2^2 =  -4 k$. 
\smallskip

From $\rho_*(B_1)=A_1$, the projection formula implies that
$$
\phi_1^*(W)\cdot B_{1} = p_1^* (W) \cdot A_1 = {k} \, ,
$$
and similarly $\phi_1^*(W)\cdot B_{2}  = p_1^* (W) \cdot A_2= {k}$.

Finally, by the Riemann-Hurwitz formula, $\omega$ is isomorphic to  
$$
B_{1} + B_{2} + \, \phi_1^*(\omega_1) = B_{1} + B_{2} + 4 \, \phi_1^*(W) \, .
$$ 
We have $\phi_1^*(W)^2 = \deg_K(<\phi_1^*(W),\phi_1^*(W)>)=0$, $B_1^2= A_1^2/2=-2k$ and  $B_2^2= A_2^2/2=-2k$, so we get 
$$
\omega^2 = B_{1}^2 + B_{2}^2 + 8 (B_{1} + B_{2}) \cdot \phi_1^*(W)
= -4k+ 16 {k} = 12 {k} \,  ,
$$

\bigskip

\subsection{\\}

The smooth curve $X$ over $Y$ is $2$-pointed since $\phi : X \to Y$ induces an isomorphism $E_i \to Y$
when $i=1,2$. The line bundle $\omega_{X/Y}+B_1+B_2$  is $\phi$-ample because, for every $\overline y \in Y(\overline K)$,  the degree of its restriction to the fiber $X_{\overline y}$ is positive. So $\phi: X\to Y$
is a stable smooth $2$-pointed curve of genus six. Let $\overline Y$ be  a regular projective model of $Y$
over $S= {\rm Spec}({\mathcal O}_K)$. By Proposition~\ref{prop19} there exists an alteration
$a_S: Y_S \to \overline Y$,  finite and \'etale over $Y$, and a stable $2$-pointed curve 
$\phi_S : X_S \to Y_S$ of genus six extending the curve $a_S^{-1}(\phi)= \phi' : X' \to Y'$.

\section{A conjecture of Parshin and Moret-Bailly}

\subsection{ \ }  We keep the notation of   §5.1. Let  $\overline K$ be an algebraic closure of $K$.   We choose an  integrable   adelic metric on  
the relative dualizing sheaf $\omega$  of $\phi : X \to Y$.
Let $F \subset \overline K$ be a finite extension of $K$
and let $P\in Y(F)$ be a rational point. Denote by  $X_{P}$ the  fiber of  the projection $\phi $ at $P$, and let $\overline\omega_P$ be the pull-back of $\overline \omega $ to $X_P$.
 We let $d$ be the absolute degree of $F$. Recall that $\delta(F)= \log(|\Delta_{F}|)$. The following is a special case of  a conjecture of Parshin \cite {Pa1} and Moret-Bailly
\cite{MB3}:

\medskip

\begin{thm}\label{thm1}
For every positive rational number $\varepsilon$, one has the inequality
$$
\overline\omega_{P}^2 \leq (36 + \varepsilon) \, \delta(F) + 0(d) \, .
$$
\end{thm}

\medskip

\noindent {\bf Remark}:
If Theorem~\ref{thm1} is valid for one choice $\overline \omega$ of an integrable adelic 
metric on $\omega$, it is also true for any other choice ${\overline \omega}'$. To see that,
 define  $\hat x \in \widehat{\rm Pic} \, (Y)$ by the 
formula
$$
\hat x= \langle {\overline \omega}',{\overline \omega}' \rangle -  \langle {\overline \omega},{\overline \omega} \rangle \, .$$
The degree of the  image $x $ of $\hat x$ in $ {\rm Pic} (Y)$
is 
$${\rm deg}_K(x)= \omega^2 -  \omega^2 =   0 \, .$$

The Deligne pairing commutes  
with base change (\cite{YZ} Thm. 4.1.3.), so we have
$$
\widehat{\deg}_F \, (u^* (\langle \overline\omega, \overline \omega \rangle))=   \overline\omega_P^2 \,  . 
$$
We have  $\omega^2 > 0$ by  Proposition~\ref{prop5}, therefore $\langle \omega, \omega \rangle$
is ample. So, by Corollary~\ref{cor1},
$$(\overline\omega_{P}')^2 \leq (1+\varepsilon) \overline\omega_{P}^2 + 0(d) \leq (36 + \varepsilon) \delta(F) + 0(d) \, .$$

\subsection{ \ }

With the notation of §5 we
 consider the line bundles $B_1$  and $B_2$,
$$B = 2B_1 + 2B_2 \, , $$
and 
$$L=\omega_{X/Y}-B$$
on $X$.
Let $\Psi \subset Y(\overline K)$ be the set of points
$\overline y$  such that $P_1$, hence $P_2$, is exceptional. 
Let $\Phi \subset Y(\overline K)$ be the set of points $\overline y$  such that, for some choice of $\cal L$,
$H^0(Y_1,P_1+\lambda) \neq 0$ or $H^0(Y_1,P_2+\lambda) \neq 0$.
If $M$ is a line bundle on $X$ (resp. $Y_1\times Y$) 
we denote by $M_{\overline y}$ the restriction of $M$ to the geometric fiber of  $\phi$ (resp. $p$) at
 $\overline y$.
We shall use the following geometric results:
\bigskip

\noindent 
\begin{prop}\label{prop6}
\begin{enumerate}
\item[i)] The set $\Psi$ is a finite set of order $4k$.
The set $\Phi$ is finite; it contains $\Psi$.
\item[ii)]
$\deg_K(L)=6$.
\item[iii)]
$$B^2=-16k \, , $$
$$\omega_{X/Y}\cdot B =8k \, ,$$
$$ L^2 = -20k \, .$$
\item[iv)]
 $h^0(B_{\overline y}) = 3$ if $\psi({\overline y})$ is exceptional,
$h^0(B_{\overline y}) = 1$ if not.
\item[v)]$h^0(L_{\overline y}) = 4$ if $\psi({\overline y})$ is exceptional,
$h^0(L_{\overline y}) = 2$ if not.
\item[vi)]
$$\deg_K (\phi_*(L)) \geq -3k \, .$$

\end{enumerate}
\end{prop}
\bigskip

\subsubsection {\\}
\noindent {\bf Proof of i).} 
Since there exist $4$  exceptional points (see Proposition~\ref{prop13} iii))    and since $\psi$ is finite of degree $k$, we conclude that $\Psi$ is a finite set of order $4k$. When $P_1$ is a Weierstrass point  we know, by Proposition~\ref{prop3}  iii), that $H^0(Y_1,P_1+\lambda) = H^0(Y_1,P_2+\lambda)  =0$. Recall from §5.2 that $\lambda= {\cal L}_{\overline y}$, $P_1=A_{1,{\overline y}}$ and $P_2= A_ {2,{\overline y}}$.          
The line bundle $A_1+A_2$ has finitely many square roots therefore, by the semi-continuity
of cohomology, $\Phi$ is finite.
Furthermore, when $P_1$ is exceptional,  $h^0(\lambda)=1$ therefore $h^0(Y_1,P_1+\lambda) \geq 1$ and $h^0(Y_1,P_2+\lambda) \geq 1$,
so $P_1\in \Phi$.

\subsubsection {\\}
\noindent {\bf Proof of ii) and iii)}

We have
$$\deg_K(L) = \deg_K(\omega_{X/Y})-4 = 10-4 = 6 \, .$$

Using Proposition~\ref{prop5} we compute
$$B^2= 4B_1^2 + 4B_2^2= -8k-8k=-16k \, ,$$
$$\omega_{X/Y}\cdot B = ({\phi_1}^*(\omega_1)+B_1+B_2)\cdot (2B_1+2B_2)  $$
$$ = 8 ({\phi_1}^*(W) \cdot (B_1+B_2)) +2 (B_1^2+B_2^2) = 16k-4k -4k= 8k           \, ,$$

and
$$L^2 = \omega_{X/Y}^2 -2\omega_{X/Y}\cdot B +B^2 = 12k-16k-16k=-20k \, .$$

\subsubsection {\\}
\noindent {\bf Proof of iv) and v).}

By Serre-duality and the Riemann-Roch theorem
$$h^0(L_{\overline y}) = h^1(B_{\overline y}) = h^0(B_{\overline y}) -4+6-1 = h^0(B_{\overline y})+1\, .$$
Since $B=\rho^*(A_1+A_2)$ we get
$$h^0(B_{\overline y}) = h^0(P_1+P_2)+h^0(\lambda) \, .$$

When $P_1$ is exceptional $h^0(\lambda)=1$ and $h^0(P_1+P_2)=2$
(Proposition~\ref{prop13} , ii) and iii)). Therefore $h^0(B_{\overline y}) =3$ and 
$$h^0(L_{\overline y}) =4 \, . \leqno (14)$$

If $P_1$ is not exceptional  $h^0(\lambda)=0$ and $h^0(P_1+P_2)=1$
(loc. cit.) so $h^0(B_{\overline y}) =1$ and 
$$h^0(L_{\overline y}) =2 \, . \leqno (15)$$

\subsubsection {\\}
\noindent {\bf Proof of vi).} 
Let $R^1\phi_*(B)_{\rm tor}$ be the torsion sub-${\mathcal O}_Y$-module of $R^1\phi_*(B)$.  The relative Serre duality theorem implies that
$$\deg_K(\phi_*(L)) =  \deg_K (R\phi_*(B)) - \deg_K( \phi_*(B))+\deg_K(R^1\phi_*(B)_{\rm tor})\, . \leqno(16)$$

By the Grothendieck-Riemann-Roch theorem 
$$
\deg_K R\phi_ {*}(B)= \deg _K(\phi_{*} ({\rm ch} (B) \, {\rm Td} (X/Y))^{(1)})    \,  ,  \leqno (17)
$$
where $x^{(1)}$ denotes the component of degree $1$ of an element  $x \in CH^{*} (Y) {\otimes}{\mathbb Q}$. We have
$$
\deg_K (\phi_{*} ({\rm ch} (B) \, {\rm Td} (X/Y))^{(1)}) = \frac{B^2}2 - \frac{\omega_{X/Y}\cdot B}{2} + \frac{ \omega_{X/Y}^2}{12}\, .  
$$
From iii) and Proposition~\ref{prop5} we know that
$B^2=-16k$, $\omega_{X/Y}\cdot B=8k$ and $\omega_{X/Y}^2 = 12k$. 
Therefore
$$\deg_K (R\phi_*(B)) = -11k \, . \leqno(18)$$

Since $B= \rho^*(A_1+A_2)$ we deduce 
from (13) that 
$$\phi_*(B) = p_*(A_1+A_2) \oplus p_*(\cal L) \, .$$
Let $j: D_2 \to Y_1 \times Y$ be the inclusion. There is a short  exact sequence of coherent ${\mathcal O}_{Y_1 \times Y}$-modules
$$
0 \to A_1 \to A_1+A_2 \to j_* \, j^* (A_1+A_2) \to 0 \, ,
$$
hence an exact sequence of ${\mathcal O}_Y$-modules
$$
0 \to p_* (A_1 ) \to p_* (A_1+A_2 ) \to j^* (A_1+A_2 )  \, \leqno (19)
$$
(note that $p \circ j$ is an isomorphism). 
Since $h^0(P_1)=h^0(P_1+P_2)=1$
when $\overline y \notin \Phi$, both
$p_*(A_1)$ and $p_*(A_1+A_2)$ are line bundles on $Y$. 
But $ j^* (A_1+A_2)$ is torsion free. So (19) implies that 
$p_*(A_1) = p_*(A_1+A_2)$. A similar argument gives $p_*(A_1)= \mathcal   O_Y$. Therefore
$$ p_*(A_1+A_2) = \mathcal   O_Y \, . \leqno (20)$$

If $\overline y \in Y(\overline K)$ is a generic point
$$p_*({\cal L})_{\overline y} = H^0(Y_1\times {\overline y}, \lambda)=0 \, ,$$
by Proposition~\ref{prop13} ii), and thus
 $$p_*({\cal L})=0 \, . \leqno(21)$$
From (20) and (21) we deduce that
$$\deg_K \phi_*(B)=0 \, . \leqno(22)$$

\subsubsection {\\}

When $T$ is a torsion ${\mathcal O}_Y$-module of finite type
$$\deg_K(T)= \sum_y {\ell}(T_y) \, ,$$
where $y$ runs over closed points of $Y$ and ${\ell}(T_y)$ is the length of $T_y$
( \cite{HA}, II, Ex. 6.12 (2)).

Let $y\in Y$ be a closed point, ${\mathcal O}_y$ its local ring, $\pi \in {\mathcal O}_y$
a uniformizer, $k(y)$ the residue field and $k(\overline y)$
an algebraic closure of $k(y)$, where $\overline y \in Y(\overline K) $.  If $M$ is an ${\mathcal O}_Y$-module we let $M[\pi]$ (resp. $M/\pi$) be the kernel
(resp. the cokernel) of the product by $\pi$ in $M_y$. 
Since $B$ is flat over $Y$ we have an exact sequence
of ${\mathcal O}_Y$-modules
$$0 \rightarrow B \xrightarrow{\pi} B\rightarrow B/\pi = B|_{X_y} \rightarrow 0 \, ,$$
where $X_y$ is the closed fiber of $\phi$ at $y$, from which we deduce a short exact sequence
of $k(y)$-vector spaces
$$ 0 \to \phi_*(B)/\pi \to H^0(X_y,B/\pi) \to R^1\phi_*(B) [\pi] \to 0 \, . \leqno (23)$$

Let $\overline y\in Y(\overline K)$ be the geometric point associated to $y$. We know from
iv) that $h^0(X_{\overline y}, B_{\overline y})= 3$ if $P_1$ is exceptional, and $1$ if not.
It follows from (23) that $\phi_*(B)/\pi$ has dimension one over $k(y)$. Therefore 
$R^1\phi_*(B)[\pi]$ has dimension $2$ over $k(y)$ when $P_1$ is exceptional, and $0$
otherwise. This implies that, when $P_1$ is exceptional,
$\ell(R^1\phi_*(B)[\pi])=2$
and, using i),
$$\deg_K(R^1\phi_*(B)_{\rm tor}) \geq 2 \times ({\rm card}(\Psi)) = 8k \, . \leqno (24)$$

From (16), (18), (22) and (24) we conclude that
$$ \deg_K(\phi_*(L)) \geq -11k+8k=-3k \, .$$
This ends the proof of Proposition~\ref{prop6}.

\bigskip

\subsection{ \\ }{\bf Proof of Theorem~\ref{thm1}.}
We use the notation of \S 5.4 . Since $a:Y'\to Y$ is finite and \'etale there exists a finite set $\cal P$
of places of $K$ such that, if $U$ is the spectrum of $\cal P$-integers in $K$,
the map $a$ extends to a finite, surjective and unramified morphism of proper $U$-schemes
$a_U: Y'_U \to Y_U$. By the Chevalley-Weil theorem, if $F\subset \overline K$ is a finite extension of $K$, any point $P\in Y(F)$ can be lifted to a point $P'\in Y'(F')$,
where $F'\subset \overline K$ is an extension of $F$ of degree at most the degree $m$ of $a$ and unramified outside $\cal P$.
Let $d=[F;\mathbb Q]$ and $d'=[F':\mathbb Q]$. It follows from \cite{Se} 1.3, Proposition 4, that
$$\frac{\delta(F')} {d'} \leq \frac {\delta(F)}{d} + {C} \, , \leqno (25)$$
where 
$$C= \sum_{p\in {\cal P}_0} \log(p) + {\rm card}({\cal P}_0)  \log(m)$$
and ${\cal P}_0$ denotes the set of residual characteristics of places in $\cal P$
(see \cite{MB3} 4.2.2)

Let ${\rm Spec}(F') \to Y$ be the map associated to the image of $P'$ 
in $Y(F')$. If $X_{P'}$ is the fiber product of $ {\rm Spec}(F')$ and $X$ over $Y$,   we let
$\overline \omega_{P'}$ be the inverse image of $\overline \omega$ on  $X_{P'}$.
We shall prove in §6.3.8 below that, for every positive rational number $\varepsilon$,
$$\overline \omega_{P'}^2 \leq (36 + \varepsilon)\delta(F')  + 0(d') \, . \leqno (26)$$
Since  $X_{P'}$ is the fiber product of  $ {\rm Spec}(F')$ and $X_P$ over $ {\rm Spec}(F)$, 
$\overline \omega_{P'}$ is the inverse image of $\overline \omega_P$ on  $X_{P'}$.
and the projection formula gives
$$\overline \omega_{P'}^2 = \frac{d'}{d} \overline \omega_{P}^2 \, . \leqno (27)$$
Theorem~\ref{thm1} follows from (25), (26) and (27).

\bigskip

\subsubsection{\\}
Let $k'$ be the degree of the composite $\psi'$ of $a: Y'\to Y$ with $\psi: Y\to Y_1$. Let ${\overline  y}' \in Y'(\overline K)$ be the image of $P'$.
If $i=1$ or $2$ denote
by $B'_i$ (resp. $A'_i$) the inverse image of $B_i$ (resp. $A_i$) on $X'$ (resp. $Y_1\times Y'$).
Let $\Psi'=a^{-1}(\Psi)$ and $\Phi'=a^{-1}(\Phi)$. Denote by $\lambda'$
the inverse image of $\lambda$ on $Y_1\times {\overline y}'$.
If $\omega'$ is the sheaf of relative differentials of $X'$ over $Y'$, we consider
the line bundles $B'= 2B'_1+2B'_2$ and
$$ L'= \omega'-B' \, $$
on $X'$
It follows from the projection formula and flat base change that Proposition~\ref{prop5} and Proposition~\ref{prop6}
remain true when $k$ (resp. $X$, resp. $Y$, resp. ${\overline y}$,  resp. $A_i$, resp. $B_i$, resp. $B$, resp. $\Psi$, resp. $\Phi$, resp. $\omega$, resp. $L$) are replaced by
$k'$ (resp. $X'$, resp. $Y'$, resp. ${\overline y}'$, resp. $A'_i$, resp. $B'_i$, resp. $B'$, resp. $\Psi'$, resp. $\Phi'$, resp. $\omega'$, resp. $L'$).

\medskip

Since $X_S$ is a $2$-pointed curve it carries two line bundles $B_{1S}$ and $B_{2S}$.
Define
$$L_{S} = \omega_{X_S/Y_S}-2B_{1S}-2B_{2S} \, .$$
The restriction of $L_{S}$  to $X'$ is the line bundle
$L'$.

\bigskip

\subsubsection{\\}

 The curve $X_{P'} $ is the fiber product of $X'$ and ${\rm Spec}(F')$ over $Y'$.
Let $v':  X_{P'} \to X'$
be the inclusion
and let $\phi'_{P'} : X_{P'} \to  {\rm Spec}(F')$,
be the projection.
 Consider the line bundle
$$L'_{P'} = v'^*(L')$$
on $X_{P'} $, and let
$$H= \phi'_{P' *}(L'_{P'}) \, $$
on ${\rm Spec}(F')$.
If   $\overline y' \notin \Phi'$ it follows from Proposition~\ref{prop6} v)
that $H$ has dimension $2$ over $F'$.

\medskip
\subsubsection{\\}

Assume  $\overline y' \notin \Phi'$. Let $s_1$ and $s_2$ be two sections in $H$. Assume that $s_1$ and $s_2$
are linearly independent. 
If $Q\in V$ was a common zero of $s_1$ and $s_2$, for every section $s\in H\otimes \overline K$
we would get $s(Q)=0$,
hence
$$h^0(L'_{P'}-Q)=h^0(L'_{P'})=2 \, ,$$
and $Q$ would be a base point of $L'_{P'}$.

By the argument of §3.5.1 there exists a divisor $D$ of degree one over $Y_1\times {\overline y'}$ such that
$2D=P'_1+P'_2$ and $D+R_i\notin \Theta$ when $i=1$ and $i=2$.
The restriction map ${\rm Pic}(Y_1\times Y')[2] \to {\rm Pic}(Y_1 \times {\overline y'})[2]$
is surjective therefore there exists $\cal L '$ on $Y_1\times Y'$
such that $2 {\cal L}' = A'_1+A'_2$ and  ${\cal L}'_{\overline y'}=  {\mathcal O}_{Y_1}(D)$.
Since ${\overline y'} \notin \Phi'$  we get $h^0(D+P'_1)=h^0(D+P'_2)=0$. From Lemma~\ref{lm15} it follows
that, if $\rho': X_{P'} \to Y_1\times {\overline y}$ is the projection defined by $D$,
${\rho'}^*(\omega_1-D) = L'_{P'}$ has no base point.
So the sections $s_1$ and $s_2$ have no common zero.

We can thus apply Proposition~\ref{prop12} to the curve $X_{P'}$ and to the sections $s_1$ and $s_2$ of $L'_{P'}$.
Choose an admissible integrable metric $\Vert \cdot \Vert_a$ on $L'$, and let $(\Vert \cdot \Vert_{a,\sup,w})$
be the associated sup-norm on $H$. 
We get the inequality
$$-{\overline  L'}^2 \leq \deg_{F'}(L'_{P'})(h_a(s_2)+ h_a(s_1) )\,  ,\leqno (28)$$
where, for $i=1$ or $2$,
$$h_a(s_i)= \sum_{w} \log (\Vert s_i \Vert_{a,\sup,w}) \, .$$

\bigskip

\subsubsection{\\}

Let $\Vert \cdot \Vert$ be the metric on $L'$ induced by the model $L_{S}$ over $S$. 
By  Lemma~\ref{lm1}, for every place $v$ of $K$, there is a constant
$C_v\geq 1$, independent of $F$, $P$, $F'$ and $P'$, such  that
$C_v=1$ for almost all $v$ and, for all $v$,
$$\Vert\cdot\Vert_{a,v} \leq  C _v \Vert \cdot\Vert_{v} \,  .$$
If $(\Vert\cdot\Vert_{\sup,w})$  is the sup-norm induced 
by $\Vert\cdot\Vert$ on $H$, and if 
$$h(s_i)= \sum_{w} \log (\Vert s_i \Vert_{\sup,w})$$
we obtain
$$h_a(s_i) \leq h(s_i) +0(d') \,  . \leqno(29)$$

\bigskip

\subsubsection{\\}

The metric  $\Vert \cdot \Vert$  on $L'$ has a model on the whole of  ${\rm Spec}({\mathcal O}_{K})$.
Therefore, by Proposition~\ref{G},  the metric  $\Vert\cdot\Vert_{\sup}$ on $H$ is adelic,
with a model on the whole of ${\rm Spec} ({\mathcal O}_{F'})$.
Assume $\overline y' \notin \Phi' $.
Then $H$ has dimension $2$ and,  by Minkowski's theorem (cf. (2)), there exist two   sections $s_1$ and $s_2$
in $H$ which are linearly independent and such that
$$h(s_1) + h(s_2)  \leq - {\widehat \deg}_{F'}(H_{\sup})+ {\delta(F') }+ A d' \, , \leqno(30)$$
where $A$ is a an absolute positive constant. 
From (28), (29), and (30), since $\deg_{F'}(L'_{P'}) =  6$, we deduce that, if  $ \overline y' \notin \Phi'  $, we have

$$-{\overline  L}'^2 \leq - 6 {\widehat \deg}_{F'} 
(H_{sup}) + 6 \delta(F') + 0(d') \, . \leqno(31)$$

\bigskip

\subsubsection{\\}

The cartesian square

$$
\xymatrix{
X_{P'} \ar[d]_{\Phi'_{P'}} \ar[r]^{v'} &X' \ar[d]^{\phi'} \\
{\rm Spec}(F')\ar[r]^-{u'} &Y' \, 
} 
$$
induces a morphism
$$\alpha : u'^*(\phi'_*(L') )\rightarrow H \, $$
(\cite{HA},  III, Remark 9.3.1). When $\overline y' \notin \Phi' $  we have
$\dim(H)=2$  so, by Grauert's theorem   (\cite{HA}, III, Cor. 12.9), $\alpha$ is an isomorphism.
Furthermore $\alpha$ is compatible  with sup-norms. Therefore $\alpha$ is an isometry. This implies that
$${\widehat \deg}_{F'}(H_{\sup})={\widehat \deg}_{F'} ( u'^*(\phi'_*(L') _{\sup})) \, .\leqno(32)$$ 

Let ${\widehat x} \in {\widehat {\rm Pic}}(Y')$ be the class of $-\phi'_*(L')_{\sup}$. The degree of  its image $x$ in  ${\rm Pic}(Y')$
satisfies
$$\deg_{F '}(x) =-\deg_{F'} (\phi'_*(L')) \leq 3k' =\frac{ \omega'^2}{4} \, ,$$
where $\omega'=\omega_{X'/Y'}$
(Proposition~\ref{prop6}  vi), and Proposition~\ref{prop5} ). Therefore, by Corollary~\ref{cor1},
$${\widehat \deg}_{F'}(u'^*({\widehat x})) \leq (1+ \varepsilon)  \frac {{\overline \omega_{P'}}^2}{4}  + 0(d') \, .\leqno(33)$$
From (31), (32), and (33) we deduce that, if $\overline y' \notin \Phi' $ ,

$$-{{\overline  L}'_{P'}}^2 \leq  (3+\varepsilon) \frac{\overline \omega_{P'}^2}{2} + 6\delta (F') + 0(d') \, . \leqno (34)$$

\bigskip

\subsubsection{\\}

Let $\widehat x \in {\rm  \widehat Pic}(Y')$ be class of the Deligne pairing
$$ \widehat x = \langle {\overline  L}' , {\overline  L}' \rangle \, .$$
Its image $x \in {\rm Pic}(Y')$ has degrre
$$\deg_{F'}(x) =L'^2= - 20 k' = -\frac 53 \omega'^2$$
(Proposition~\ref{prop6} iii) and Proposition~\ref{prop5}). Therefore, by Corollary~\ref{cor1}, 
$${\widehat \deg}_{F'} ( u'^*({\widehat x}))  = {{\overline  L}'_{P'}}^2 \leq (-\frac 53+ \varepsilon) \, {\overline \omega}_{P'}^2 + 0(d') \, . \leqno(35)$$
From (34) and (35) we deduce that, if $\overline y' \notin \Phi'  $,
$$\frac 53 \overline {\omega}_{P'}^2  \leq  (3+\varepsilon) \frac{\overline {\omega}_{P'}^2 }{2} + 6 
\delta(F')+0(d') \, ,$$
i.e., with another $\varepsilon$,
$$\overline{ \omega}_{P'}^2 \leq ( 36 +\varepsilon) \delta(F') + 0(d') \, . \leqno(36) $$

\bigskip

\subsubsection{\\}
By the projection formula, the quantity $c(\overline y') =  \overline\omega_{P'}^2/{d'} $
does not depend on $P'$ but only on its image $\overline y'$ in $ Y'(\overline K)$.
Let
$C$ be the supremum of  $|c(\overline y')|$ when $\overline y'$
runs over the finite set $ \Phi' $.
By (36)  we get, for all $P'$,
$$
\overline\omega_{P'}^2
 \leq ({36}+ \varepsilon) \delta (F') + 0(d')+  C d'  = ({36} + \varepsilon) \delta (F') + 0(d')\, .  $$
This ends the proof of (26), hence the proof  of Theorem 1.

\bigskip

\section{Consequences}\label{sec3}
\subsection{ \\ } {\bf  Effective Mordell for curves of genus two.} 
\subsubsection{\\}
Let $K_0$ be a number field and ${\mathcal C}$ a geometrically
 connected, smooth and projective curve of genus $2$ over $K_0$. Let $\overline K_0$ be an algebraic closure of $K_0$ and let
$$
h : {\mathcal C}(\overline K_0) \to {\mathbb R}
$$
be a logarithmic height associated to the canonical line bundle over ${\mathcal C}$.

\begin{thm}\label{thm2} For every  rational number $\varepsilon > 0$ there exists a constant $\beta_1(\varepsilon) > 0$ such that, 
for any finite extension $F$ of $K_0$ and $P \in  {\mathcal C}(F)$,
$$
h(P) \leq (12 +\varepsilon)  \frac{\delta (F)}{[F:{\mathbb Q}]} + \beta_1(\varepsilon) \, .
$$

\end{thm}

\bigskip

\subsubsection{ \ }  {\bf  Proof.}  We use the notation of  §5.1  and     §5.3. 
 Let $\pi: Y \to  {\mathcal C}$
be the composite of the morphisms $ Y \to Y_0$ and $Y_0 \to  {\mathcal C}$.
The morphism $\pi$ is \'etale, finite and surjective. Let $M$ be the canonical line bundle on $ {\mathcal C}$.
Choose  integrable adelic metrics on $M$ and $\omega$,  and let $\overline L = \pi^* (\overline M)$.
Define
$$
\hat x= \overline L -   \frac{1}{3} \langle \overline \omega , \overline \omega \rangle 
$$
in $\widehat{\rm Pic} \, (Y)_{\mathbb Q}$. The image $x$ of $\hat x$ in ${\rm Pic} \, (Y)_{\mathbb Q}$ is such that
$$
\deg _K(x) = \deg_K(L) - \frac{1}{3} \,  \omega^2 \, . \leqno(37)
$$
We deduce from Proposition~\ref{prop5} that
$$
\omega^2 =  6 \deg \, (Y/Y_0) \, .  \leqno(38)
$$
Therefore, since $\deg_K(L) = \deg_{K _0} (M)   \deg \, (Y/Y_0)$ and $\deg_{K _0} (M) = 2 $, we get
$$
\deg_K (x) = 0 \,  .$$
If $F' \supset K$ is a finite extension and, if $P' \in Y(F')$ defines
$
u : {\rm Spec}(F') \to Y 
$, 
Corollary~\ref{cor1} implies that
$$
\widehat{\deg}_{F' }\, (u^* (\hat x)) \leq \varepsilon {\overline \omega}_{F'}^2 + 0(d') \, , \leqno (39)
$$
where  
$
d' = [F':{\mathbb Q}]  
$.
In other words,
$$\widehat{\deg}_{F'} \, (u^* (\overline L)) \leq  
\frac {1+\varepsilon}{3} \,{\overline \omega}_{F'}^2 + 0(d') \, .  \leqno (40)$$
We know from Theorem~\ref{thm1}  that
$$
\overline \omega_{F'}^2 \leq (36 +\varepsilon) \, \delta (F') + 0(d')  \, , \leqno(41)
$$
so (40) implies that
$$
\widehat{\deg}_{F'} \, (u^* (\overline L)) \leq   \left( 12 + \varepsilon \right) \delta(F') + 0(d') \, .
 \leqno (42)
$$
If $h_L$ is a  logarithmic height associated to $L$, we have
$$
 d' \, h_L (P')  \leq \widehat{\deg}_{F'} \, (u^* (\overline L)) + 0(d') \, , \leqno(43)$$
so, if we divide (42) and (43) by $d'$, we get
$$
h_L (P') \leq  \left(12 + \varepsilon \right) 
\left(\frac{\delta(F')}{d'} + 0(1)\right) \, , \leqno (44)$$
where $0(1)$  is a constant independent of $F'$  and $P'$, but depending on $\varepsilon$.

As we saw in §6.3,
the point $P \in  {\mathcal C}(F)$ can be lifted to a point $P'\in Y(F')$,
where $F'$ is a finite extension of $F$ such that, if $d= [F:\mathbb Q]$,
$$ \frac  {\delta(F')} {d'} \leq  \frac  {\delta(F)} {d} +0(1) \, . \leqno(45)$$
Furthermore, if $ h$ is  a height associated to $M $, the projection formula implies that
$$
h_L (P') = h (P) + 0(1) \, . \leqno(46) $$
We conclude from (44), (45) and  (46)  that
$$
h (P) \leq \left( 12 +
 \varepsilon \right) \frac{\delta (F)}{d} + \beta_1(\varepsilon) \, ,
$$ 
and Theorem~\ref{thm2} is proven.

\bigskip

\subsubsection {\bf Remark.} It follows from Vojta's  Conjecture 5.2.6 in \cite{V} that
$$
h(P) \leq (1 +\varepsilon) \, \frac{\delta (F)}{[F:{\mathbb Q}]} + \beta_1(\varepsilon) \,  .
$$

\bigskip

\subsection{\\}  {\bf The conjecture of the discriminant.} 

\bigskip

Let $E$ be an elliptic curve over ${\mathbb Q}$, with semi-stable reduction over ${\mathbb Z}$, let $\Delta (E)$ be the minimal discriminant of $E$, and let $N(E)$ be its conductor.

\medskip

\begin{thm}\label{thmDI}
For every rational number $\varepsilon > 0$ there exists a constant $\beta_2(\varepsilon) > 0$ such that
$$
\log \vert \Delta (E) \vert \leq ( 144  \,  ( 1+\log( 24) )+  \varepsilon) \log N(E) + \beta_2(\varepsilon) \, .
$$
\end{thm}

\medskip

\noindent{\bf Proof.} 
To prove Theorem~\ref{thmDI}, we proceed as in \cite{MB3}, \S6. Consider the modular curve $C = X_0 (23)$. Its genus is 2. The degree of the $j$ homomorphism
$$
j : X_0 (23) \to {\mathbb P}^1 \, ,
$$
{\it i.e.} the index of $\Gamma_0 (23)$ in $SL_2 ({\mathbb Z})$, is equal to $n=24$ (\cite{Ogg}, p. IV.4.). The degree of $L = j^* \, {\mathcal O} (1)$ is equal to $n \deg ({\mathcal O} (1)) = n$. If $h$ is the na{\"\i}ve height of ${\mathbb P}^1$ and $h_1 = h \circ j$, the height $h_1$ is a logarithmic height associated to $L$.
 Theorem~\ref{thmDI} follows from these facts, from Theorem~\ref{thm2} applied to $C$, and from \cite{MB3}, Theorem~6.4.

\bigskip

 \subsection{ \\} {\bf  A  weak uniform $abc$ theorem.} 

\subsubsection{\\}

Let  $r$ be the  degree of $K$ over ${\mathbb Q}$. For any finite place $v$ of $K$ and $x \in K^*$, we denote by $\Vert x \Vert_v$ the norm of $x$ at $v$, normalized by the equality
$$
\Vert p \Vert_v = p^{-1/r}$$
if the prime $p$ divides $v$, and we denote by $Nv$ the cardinality of the residue field. 
When $v$ is a real place we let
$$
\Vert x \Vert_v = |x|^{-1/r}
$$
and when $v$ is complex we let
$$
\Vert x \Vert_v = |x|^{-2/r} \, .
$$

 Given three elements $a,b,c$  in $K^*$ we define
$$
H(a,b,c) = \prod_v \max (\Vert a \Vert_v , \Vert b \Vert_v , \Vert c \Vert_v) \leqno(47)
$$
and 
$$
N(a,b,c) = \prod_{v \in I} Nv^{1/r} \leqno(48)$$
where $I$ is the set of finite places of $K$ such that
$$
\max (\Vert a \Vert_v , \Vert b \Vert_v , \Vert c \Vert_v) > \min (\Vert a \Vert_v , \Vert b \Vert_v , \Vert c \Vert_v) \, .
$$

\begin{thm}\label{thm3}
For every rational number $\varepsilon > 0$ there exists a positive constant $\beta_3 ( \varepsilon)> 0$ independent of $K$
 such that, if $a+b = c$,
$$
\log H(a,b,c) \leq  ( 24   + \varepsilon) \log N(a,b,c) + ( 30 +\varepsilon) \frac{ \log (|\Delta_K|)}{[K:{\mathbb Q}]}+ \beta_3 ( \varepsilon)\, .
$$

\end{thm}

\bigskip

\subsubsection{ \\ }  {\bf Proof.} 
To prove Theorem~3, we proceed as in \cite{MB3}, \S7. Let $\mathcal C$ be the smooth
 projective curve over $K_0={\mathbb Q}$ which is the complete model of the affine curve with equation
$$
y^2 + y = x^5 \, .
$$
This curve $\mathcal C$ is geometrically connected and  it has genus 2.

 The coordinate $y$ defines a finite surjective map
$$
\pi : \mathcal C \to {\mathbb P}^1
$$
of degree 5, totally ramified at $0$, $-1$ and $\infty$.  Le $L= \pi^* \, {{\mathcal O} (1)}$.
The projection formula implies
$$
\deg_{K_0}(L) = \deg_{K_0} (\pi^* \, {{\mathcal O} (1)}) =  \deg_{K_0} (\pi_* \pi^* \, {{\mathcal O} (1)})
= 5 \deg_{K_0} ({\mathcal O} (1)) = 5 \, .$$

Let $y = \frac{b}{a} \in {\mathbb P}^1 (K)$, $m = y^2 + y = bc / a^2$, 
and $F = K(\sqrt[5]{m})$. Consider the point $P = (\sqrt[5]{m} , y) \in {\mathcal C}(F)$. If $h$ is the na{\"\i}ve normalized height on ${\mathbb P}^1(K)$ and if $h_{ L}$ is a logarithmic  height attached to  $ L$, we have
$$
h(y) \leq h_{L} (P) +0(1) \, ,
$$
hence, by Theorem~\ref{thm2}, since $\frac {\deg_{K_0}(L) } {2}\times 12 =30$,
$$h(y) \leq (30+ \varepsilon)\, \frac{\log (|\Delta_F|)}{[F : {\mathbb Q}]} + 
\beta_1 ( \varepsilon) \, \leqno  (49)$$
for some constant $\beta_1 ( \varepsilon) >0$. We have
$$
 \log H (a,b,c)  \leq h(y) + 0(1) \, .
\leqno(50)
$$
If $\Delta_{F/K}$ is  the norm of the discriminant of $F$ over $K$, we get
$$
|\Delta_F| = |\Delta_{K}|^n \, \Delta_{F/K} \, ,
$$
where $n$ is the degree of $F$ over $K$. From \cite{Se} (5), we deduce that
$$
\log \, (\Delta_{F/K}) \leq [F:{\mathbb Q}] \left( 1-\frac{1}{n}  \right) \log \, N(a,b,c) +   [F:{\mathbb Q}] \,\log(5) \,  . 
$$
Since $n=1$ or $5$ we get
$$
\frac{\log (|\Delta_F|)}{[F:{\mathbb Q}]} \leq \frac{4}{5} \log N(a,b,c) +  \frac{\log (|\Delta_K|)}
{[K :  {\mathbb Q}]} + \log(5) \, .
\leqno(51)  $$
From   (49), (50) and  (51) we deduce that, for every rational number $\varepsilon > 0$, there exists a constant $\beta_3 ( \varepsilon) > 0$ such that
$$
\log H(a,b,c) \leq  ( 24   + \varepsilon) \log N(a,b,c) + ( 30+\varepsilon) 
\frac{ \log (|\Delta_K|)}{[K:{\mathbb Q}]}+ \beta_3 ( \varepsilon)\, .
$$

\hfill q.e.d.

\bigskip

\bigskip

\subsection{ \ }{\bf  No Siegel zero.}

Let $-d < 0$ be a fundamental discriminant. Consider the quadratic character $\chi_d := \left( \frac{-d}{\cdot} \right)$ of conductor $d$, and the associated Dirichlet $L$-function $L(s,\chi_d)$.

\bigskip

\begin{thm}\label{thm77}
 There exists a real constant $c > 0$ such that $L (s,\chi_d)$ has no real zero in the interval
$$
1 - \frac{c}{\log (d)} < s \leq 1 \, .
$$
\end{thm}

\bigskip

\noindent {\bf Proof.} Theorem 6 follows from the uniform $abc$ conjecture, as was shown by Granville and Stark  \cite{G-St}.  Elkies  explained  that their argument remains valid with a different exponent
than $1+\varepsilon$.

If $A$ and $B$ are positive real functions we write $A<<B$ (resp. $A\asymp B$) when there exists a positive constant  $c$
such that $A \leq cB$ (resp. when $A<<B$ and $B<<A$). 

We let
$$d_K = |\Delta_K|^{\frac{1}{[K:\mathbb Q]}} \, .$$  Theorem~\ref{thm3} implies that, if $a,b,c$ are in $K$ and
$a+b=c$, then
$$H(a,b,c) << N(a,b,c)^{25} d_K^{31} \, .$$

Let $E$ be an ellliptic curve with complex multiplication by the ring of integers of $k=\mathbb Q(\sqrt{-d})$ and equation
$$y^2=4x(x-1)(x-\lambda) \, .$$
Let $j$ be the $j$-invariant of $E$ and $L=k(j,\lambda)$.
It is proved in \cite{G-St} \S 4.1 that one can choose $\lambda$ in  such a way
that $N(\lambda,1-\lambda,1) \asymp 1$, and that $d_L<< \sqrt d$. Therefore

$$H(\lambda, 1-\lambda,1) << d^{16} \, .$$

Let $f: \mathbb P^2 \to \mathbb P ^1$ be the rational map such that
$$ f(x_0, x_1, x_2) = (256(x_0^2+x_1x_2)^3,(x_0x_1x_2)^2) \, .$$
The map $f$ has degree $6$ so, when $ f(P)$ is well defined, the Weil height of $f(P)$
verifies
$$ H(f(P)) << H(P)^6 \, $$
(\cite{HS} Theorem B.2.5., (a)). If $x_0=\lambda$, $ x_1 = 1-\lambda $, and $x_2=1$, the map $f$ is well
defined at $P= (x_0,x_1,x_2)$. By \cite{G-St}$ (12)$ one gets
$$H(j,1) = H(256(\lambda^2-\lambda+1)^3, (\lambda^2-\lambda)^2) << H(\lambda,1-\lambda, 1)^{6} _, .$$
Therefore
$$H(j,1) << d^{96} \, .$$

Let $h$ be the class number of $\mathbb Q(\sqrt {-d })$ and
$$h'= {h}/{\sum_{(a,b,c)} \frac{1}{a}}  \, ,$$
where $(a,b,c)$ runs over all reduced  integral quadratic forms $ax^2+bxy+cy^2$ of discriminant $-d = b^2-4ac$.
According to \cite{G-St}$ (7)$
$$H(j,1) \asymp \exp(\frac{\pi \sqrt d}{h'} )\, .$$ 
Therefore, when $d$ is large enough, 
$$h' \geq  \Gamma \frac {\sqrt d}{\log (d)} \, ,  \leqno(52)$$
with $\Gamma = \frac{\pi}{96}$ .

Let $c_2$, $c_6$ and $\sigma_1$ be the positive constants defined in \cite{MA} \S 1, and 
$\gamma = \frac {c_6}{c_2} \Gamma$. It is shown in op.cit. \S3, that, when (52) is satisfied, there exists some positive integer $d_0$ such that, if
$d > d_0$ and 
$$\max(\sigma_1, 1- \frac{\gamma}{\log(d)}) < s \leq 1 \, ,$$
then $L(s,\chi_d)$ is positive. For every $d \leq d_0 ,  d\geq 2$, choose $\sigma(d) < 1$ such that, if $\sigma(d) < s \leq 1$,
$L(s,\chi_d)$ is positive (such a $\sigma(d)$ exists since $L(1,\chi_d)$ is positive). Let
$$\sigma_2 = \sup_{2 \leq d\leq d_0} \sigma(d) \, ,$$
$$ \gamma_1 = (1-\sigma_1) \log(2) \, ,$$
$$\gamma_2 = (1-\sigma_2) \log(2) \, ,$$
and
$$c = \inf(\gamma_1,\gamma_2, \gamma) \, .$$
When $-d$ is a fundamental discriminant and
$$1-\frac{c}{\log(d)} < s \leq 1 \, ,$$
$L(s,\chi_d)$ is positive. This proves Theorem~\ref{thm77} .

\end{document}